\documentclass[10pt]{article}

\def\dofigures{11}

\usepackage{undertilde}

\usepackage{calc}
\usepackage{color}
\usepackage{amsmath}
\usepackage{amssymb}
\usepackage{mathrsfs}
\usepackage{amsthm}
\usepackage[sans]{dsfont}

\usepackage{textcomp}

\usepackage{bbding}
\usepackage{perpage}

\usepackage{mathtools}

\usepackage{fancyhdr}
\usepackage{supertabular}

\usepackage{verbatim}

\usepackage{alphalph}

\theoremstyle{plain}
\newtheorem{theorem}{Theorem}

\newtheorem{lemma}{Lemma}

\theoremstyle{definition}

\theoremstyle{remark}
\newtheorem{remark}[lemma]{Remark}

\usepackage{perpage}

\MakePerPage{equation} \renewcommand{\theequation}{\theperpage\alphalph{\value{equation}}\relax} %
\MakePerPage{theorem}  %
\MakePerPage{lemma}  %
\MakePerPage{proposition}  %
\MakePerPage{remark}  %
\MakePerPage{definition}  %
\MakePerPage{example}  %
\MakePerPage{corollary}  %

\newtagform{brackets}{\textlbrackdbl}{\textrbrackdbl}
\usetagform{brackets}

\newcommand{\myref}[1]{\ref{#1}\relax}
\newcommand{\myeqref}[1]{\textlbrackdbl\ref{#1}\textrbrackdbl\relax}
\newcommand{\mylabel}[1]{\label{#1}\relax}

\newcommand{\myeqlab}[1]{\tag{\theequation}\stepcounter{equation}\label{#1}}

\newif\ifpdf
\ifx\pldfoutput\undefined \pdffalse
  \else \pdfoutput=1 \pdftrue \fi

\ifpdf
  \usepackage[pdftex]{graphicx}
  \DeclareGraphicsExtensions{.pdf,.png,.mps}
\else
  \usepackage[dvips]{graphicx}
  \DeclareGraphicsExtensions{.ps,.eps}
\fi

\definecolor{orange}{rgb}{1,0.5,0}
\definecolor{brown}{rgb}{0.5,0.25,0}
\definecolor{violet}{rgb}{0.25,0.0,0.65}
\definecolor{darkblue}{rgb}{0,0,0.5}
\definecolor{darkgreen}{rgb}{0,0.5,0}
\definecolor{darkred}{rgb}{0.65,0,0}
\definecolor{lightgrey}{rgb}{0.6,0.6,0.6}
\definecolor{lightcyan}{rgb}{0.3,1,1}
\definecolor{medcyan}{rgb}{0,0.7,0.7}
\definecolor{midgreen}{rgb}{0,0.8,0}

\definecolor{result}{rgb}{0.5,0,0}

\def\defeq{:=}
\def\eqdef{=:}
\def\eq{=}

\newcommand{\Det}{\operatorname{det}}
\newcommand{\adjt}{\operatorname{adj}}
\newcommand{\impl}{\Rightarrow}
\newcommand{\trace}{\operatorname{tr}}

\newcommand{\topref}[2]{\overset{\text{\myeqref{#1}}}{#2}}
\newcommand{\toprefb}[3]{\overset{\text{\myeqref{#1}}}{\underset{\text{\myeqref{#2}}}{#3}}}
\newcommand{\toprefc}[4]{\overset{\text{\myeqref{#1},\myeqref{#2}}}{\underset{\text{\myeqref{#3}}}{#4}}}
\newcommand{\toprefd}[5]{\overset{\text{\myeqref{#1},\myeqref{#2}}}{\underset{\text{\myeqref{#3},\myeqref{#4}}}{#5}}}

\newcommand{\supeq}[2]{\mathord{\overbrace{\mathop{#1}}^{#2}}}
\newcommand{\subeq}[2]{\mathord{\underbrace{\mathop{#1}}_{#2}}}

\newcommand{\sign}{\operatorname{sgn}}

\newcommand{\norm}[1]{\left\|#1\right\|}
\newcommand{\norma}[2]{\norm{#1}_{#2}}
\newcommand{\anorma}[1]{\norma\cdot{#1}}

\newcommand{\osc}{\operatorname{osc}}

\newcommand{\csep}{\quad,\quad}

\newcommand{\spC}{\mathcal{C}}
\newcommand{\Cinf}{\Ck\infty}
\newcommand{\Ck}[1]{\spC^{#1}}

\newcommand{\Cone}{\Ck1}

\newcommand{\nabt}{\nabla^T}

\newcommand{\Dfun}{\mathcal{D}} 
\newcommand{\Ddist}{\Dfun'} 

\newcommand{\CC}{\mathds{C}}
\newcommand{\R}{\mathds{R}}
\newcommand{\T}{\mathds{T}}
\newcommand{\Rplus}{\R_+}
\newcommand{\K}{\mathds{K}}

\newcommand{\N}{\mathds{N}}

\newcommand{\Z}{\mathds{Z}}

\newcommand{\boi}[2]{{]#1,#2[}}

\newcommand{\roi}[2]{{[#1,#2[}}
\newcommand{\cli}[2]{{[#1,#2]}}
\newcommand{\pvint}{\text{p.v.}\int}
\newcommand{\defm}[1]{\emph{#1}}
\newcommand{\tv}{t}

\newcommand{\vorf}{w}

\newcommand{\mex}{\mx{2}}

\newcommand{\BV}{\operatorname{BV}}
\newcommand{\loc}{{\operatorname{loc}}}
\newcommand{\Ms}{\mathcal{M}}
\newcommand{\Wi}{\mathcal A} 
\newcommand{\fdomC}{\fdom_\CC}
\newcommand{\stfdomC}{\stfdom_\CC}
\newcommand{\fdom}{\mathcal{F}}
\newcommand{\brad}{\epsilon_\stf}
\newcommand{\frad}{\epsilon_\Vort}
\newcommand{\vrad}{\epsilon_{\vrdini}}

\newcommand{\Vortdom}{\mathcal{W}}

\newcommand{\cvz}{{\ssac{\vec z}}}

\newcommand{\mv}{\mu}
\newcommand{\id}{\operatorname{id}}
\newcommand{\eps}{\epsilon}
\newcommand{\supp}{\operatorname{supp}}
\newcommand{\xv}{x} 
\newcommand{\Nv}{N}
\newcommand{\nv}{n}
\newcommand{\Gam}{\Gamma} 
\newcommand{\gam}{\gamma} 

\newcommand{\Vort}{\Omega}

\newcommand{\pct}{\partial_{\ct}}

\newcommand{\eqv}{\Leftrightarrow}
\newcommand{\ta}{\beta}
\newcommand{\pola}{\theta}
\newcommand{\polb}{\vartheta}
\newcommand{\uan}{\phi} 
\newcommand{\Uan}{\varphi}

\newcommand{\rad}{r}
\newcommand{\lrv}{a}

\newcommand{\pp}{\partial_\pp}
\newcommand{\pta}{\partial_\ta}
\newcommand{\pu}{\partial_\uan}

\newcommand{\Tav}{{\vec a}}
\newcommand{\Pav}{{\vec b}}
\newcommand{\fff}{f}

\newcommand{\setdiff}{\backslash}

\newcommand{\leb}{\ell}
\newcommand{\leba}[1]{\leb^{#1}}
\newcommand{\lone}{\leba1}

\newcommand{\Cuv}{\spC_{u0}}
\newcommand{\Cuvk}[1]{\Cuv^{#1}} 

\newcommand{\Cub}{\spC_{ub}}
\newcommand{\Cubk}[1]{\Cub^{#1}}

\newcommand{\Cbk}[1]{\Cb^{#1}}

\newcommand{\Cb}{\spC_b}
\newcommand{\asy}[1]{\langle #1\rangle}

\newcommand{\embed}{\hookrightarrow}
\newcommand{\peta}{\partial_\etav}

\newcommand{\pU}{\partial_\Uan}
\newcommand{\pa}{\partial_\pola}
\newcommand{\pT}{\partial_\lrv}

\newcommand{\pt}{\partial_t}

\newcommand{\tensor}{\otimes}
\newcommand{\vv}{\vec v}

\newcommand{\stfdom}{\Psi}

\newcommand{\stfnbh}{\ballfrc{\stfdom}{\brad}{\sstfz}}
\newcommand{\Vortnbh}{\ballfrc{\Vortdom}{\frad}{\sVortz}}

\newcommand{\ballsym}{B}
\newcommand{\ballr}[1]{\ballsym_{#1}}
\newcommand{\ballfr}[2]{\ballr{#2}^{(#1)}}
\newcommand{\ballfrc}[3]{{\ballfr{#1}{#2}{(#3)}}}
\newcommand{\ballrc}[2]{\ballsym_{#1}(#2)}

\newcommand{\cballsym}{\overline\ballsym}

\newcommand{\cballrc}[2]{\cballsym_{#1}(#2)}

\newcommand{\ssac}{\check}
\newcommand{\scac}{\overline} 
\newcommand{\liac}[1]{\delta(#1)}
\newcommand{\liacs}[1]{\delta #1}

\newcommand{\zeroaccl}[1]{(#1)_0}
\newcommand{\zeroacc}[1]{#1_0}

\newcommand{\modac}{\scac}
\newcommand{\mdta}{\modac\partial_\ta}
\newcommand{\mdU}{\modac\partial_\Uan}
\newcommand{\mdu}{\pu} 
\newcommand{\mda}{\modac\partial_\pola}
\newcommand{\mdT}{\modac\partial_\lrv}

\newcommand{\cstf}{\ssac\stf}
\newcommand{\crad}{\ssac r}
\newcommand{\Vortz}{\zeroacc\Vort}
\newcommand{\sVortz}{\zeroacc\sVort}
\newcommand{\cvv}{\ssac\vv}
\newcommand{\cvort}{\ssac\vort}

\newcommand{\cvorf}{\ssac\vorf}
\newcommand{\svort}{\scac\vort}
\newcommand{\svortz}{\zeroacc{\scac\vort}}
\newcommand{\sstf}{\scac\stf}

\newcommand{\sVort}{\Vort}

\newcommand{\srad}{\scac\rad}
\newcommand{\sradz}{\zeroacc{\scac\rad}}
\newcommand{\svorf}{\scac\vorf}
\newcommand{\sstfz}{{\zeroacc\sstf}}
\newcommand{\lsstf}{\delta\sstf}
\newcommand{\lsstff}{\lsstf^\wedge}
\newcommand{\ffff}{\fff^\wedge}

\newcommand{\vrd}{h}

\newcommand{\vrdinic}{\underline\vrdini}
\newcommand{\vrdiniex}{\sigma}
\newcommand{\vrdini}{\mathring\vort}
\newcommand{\vrdinil}{\liacs\vrdini}
\newcommand{\vrdiniz}{\zeroacc\vrdini}

\newcommand{\skw}{J}

\newcommand{\va}{\vec a}
\newcommand{\vb}{\vec b}
\newcommand{\vort}{\omega}
\newcommand{\Vortl}{\dot\Vort}
\newcommand{\sVortl}{\dot\sVort}
\newcommand{\stf}{\psi}
\newcommand{\stfz}{\zeroacc\stf}

\newcommand{\linv}{\ell}

\newcommand{\pb}{\partial_\polb}

\renewcommand{\mex}{m}

\newcommand{\exs}{s}

\newcommand{\Keps}{\lops^{-1}}
\newcommand{\pole}{z}
\newcommand{\polo}{\moacs z}

\newcommand{\exso}{\moacs\exs}
\newcommand{\lopa}[2]{((\etao-\moacs{#1})\peta-\moacs{#2})} 
\newcommand{\lopad}[2]{(\peta(\etao-\moacs{#1})-(\moacs{#2}+\moacs1))} 
\newcommand{\lops}{\lopa\pole\exs}
\newcommand{\lopsd}{\lopad\pole\exs}

\newcommand{\Lop}{T} 
\newcommand{\Lopn}{T_\nv} 
\newcommand{\Eop}{E}
\newcommand{\Rop}{R}
\newcommand{\Sop}{S}

\newcommand{\ct}{{\ssac t}}
\newcommand{\vx}{{\vec x}}
\newcommand{\cvx}{{\ssac\vx}}
\newcommand{\vy}{{\vec y}}
\newcommand{\cvy}{{\ssac\vy}}
\newcommand{\cx}{\ssac x}
\newcommand{\cy}{\ssac y}
\newcommand{\vvxx}{\vv^{(\vx)}}

\newcommand{\vw}{\vec w}

\newcommand{\tao}{\acute\ta}

\newcommand{\LoneR}{\Lone(\R)}
\newcommand{\Mh}{{\LoneR+\CC\delta}}

\newcommand{\loner}{\ell^1_\Nv}

\newcommand{\Leb}{\mathcal L}
\newcommand{\Leba}[1]{\Leb^{#1}}

\newcommand{\Lone}{\Leba1}

\newcommand{\moacs}[1]{\acute#1} 
\newcommand{\moac}[1]{M(#1)}
\newcommand{\etav}{p}
\newcommand{\xiv}{\xi}
\newcommand{\etao}{\moacs{\etav}}
\newcommand{\nvo}{\moacs{\nv}}

\newcommand{\cf}{\mathds{1}}

\newcommand{\Ao}{Q}
\newcommand{\Bo}{P}

\newcommand{\normm}[1]{\norma{#1}{\Mh}}
\newcommand{\normM}[1]{\norma{#1}{[\Mh]}}

\newcommand{\mepn}{\mex_\pm}
\newcommand{\mep}{\mex_+}
\newcommand{\men}{\mex_-}

\newcommand{\mmepn}{\moacs\mex_\pm}
\newcommand{\mmep}{\moacs\mex^{}_+}
\newcommand{\mmen}{\moacs\mex^{}_-}

\newcommand{\Aopn}{\Ao+\mmepn}
\newcommand{\Aop}{\Ao+\mmep}
\newcommand{\Aon}{\Ao+\mmen}

\newcommand{\vortterm}{W}

\newcommand{\svortterm}{\scac W}
\newcommand{\gU}{g^{(\Uan)}}
\newcommand{\gu}{g^{(\uan)}}
\newcommand{\sgU}{\scac{g^{(\Uan)}}}
\newcommand{\sgu}{\scac{g^{(\uan)}}}

\newcommand{\Fmap}{F}

\newcommand{\ndim}{2}

\newcommand{\prs}{\pi}

\newcommand{\gr}{r}
\newcommand{\gs}{s}
\newcommand{\gt}{t}
\newcommand{\mgr}{\moacs\gr}
\newcommand{\mgs}{\moacs\gs}
\newcommand{\mgt}{\moacs\gt}

\newcommand{\qv}{{q}}

\begin{document}

\let\textlabel=\label

\title{Algebraic spiral solutions of 2d incompressible Euler\footnote{This material is based upon work partially supported by the
National Science Foundation under Grant No.\ NSF DMS-1054115 and by a Sloan Foundation Research Fellowship}}
\author{Volker Elling}

\maketitle

\begin{abstract}
    We consider self-similar solutions of the 2d incompressible Euler equations.
    We construct a class of solutions with vorticity forming algebraic spirals near the origin, 
    in analogy to vortex sheets rolling up into algebraic spirals. 
\end{abstract}

\section{Introduction}

\subsection{Motivation}
 
We consider the two-dimensional incompressible Euler equations: 
\begin{alignat*}{5}
    \vv_t + \nabla_\vx\cdot(\vv\tensor\vv) + \nabla\prs &= 0 \quad,\quad \nabla\cdot\vv = 0 
    \myeqlab{eq:euler}
\end{alignat*}
The divergence constraint $\nabla\cdot\vv= 0$ implies 
\begin{alignat*}{5} \vv=\nabla^\perp\stf \myeqlab{eq:stf} \end{alignat*}
for a \defm{stream function} $\stf$.
Assuming sufficient regularity we may take curl $\nabla\times$ to obtain
the \defm{vorticity}  
\begin{alignat*}{5} \vort=\nabla\times\vv=\Delta\stf \myeqlab{eq:vort-def} \end{alignat*}
which satisfies 
\begin{alignat*}{5}
    0 = \vort_t + \nabla_\vx\cdot(\vort\vv) \myeqlab{eq:vort-div} 
    \overset{\nabla\cdot\vv= 0}{=} 
    \vort_t + \vv\cdot\nabla\vort &= 0 
\end{alignat*}
In cases where the vorticity is supported on curves called \defm{vortex sheets}, 
the \defm{Birkhoff-Rott equation} is used: $W= v^x-iv^y$, the complex velocity, 
is a holomorphic function of $Z= x+iy$ in regions of zero vorticity and satisfies
\begin{alignat*}{5}
    \partial_t Z(\Gamma,t) &=  \Big(\frac1{2\pi i}\pvint\frac{d\Gamma'}{Z(\Gam,t)-Z(\Gam',t)}\Big)^*
    \myeqlab{eq:br}
\end{alignat*}
where $*$ is complex conjugation and $\Gamma$ the \defm{circulation}. 

\if\dofigures%
\begin{figure}[h]
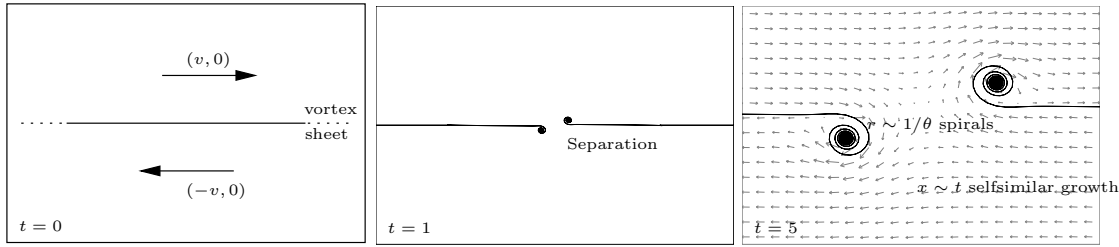

\input{pullin0.pstex_t}
\input{pullin1.pstex_t}
\input{pullin2.pstex_t}
\caption{Pullin's self-similar ($x\sim t$) separating vortex sheet}
\mylabel{fig:pullin}
\end{figure}
\fi%

Of particular interest both in applications and theory are \defm{self-similar} solutions of the Euler equations,
in particular vortex sheets:
\begin{alignat*}{5}
    Z(\Gam,t) &=  t^\mv z(\gam) \quad,\quad \Gam= t^{2\mv-1}\gam
\end{alignat*}
The shape of the sheet is the same at all times after dilation by a factor $t^\mv$. 
This ansatz yields the \defm{self-similar Birkhoff-Rott equation}
\begin{alignat*}{5}
    (1-2\mv)\gam z_\gam(\gam) + \mv z(\gam) &=  \Big(\frac1{2\pi i}\pvint\frac{d\Gamma'}{z(\gam,t)-z(\gam',t)}\Big)^*
\end{alignat*}
Prototypical solutions are vortex spirals, such as the well-known Kaden spiral \cite{kaden-1931} and Pullin's separated spirals
(Figure \ref{fig:pullin}, \cite{pullin-nuq}). Vortex sheets and spirals are observed at trailing edges of aircraft wings
and at flow past sharp corners. 
To guess a formula for such spirals, Kaden (see also \cite{rott-1956}) 
approximated the almost circular turns
near the spiral center by a single point vortex, which corresponds mathematically to 
\begin{alignat*}{5}
    (1-2\mv)\gam z_\gam(\gam) + \mv z(\gam) &=  \Big(\frac1{2\pi i}\frac{\gam}{z(\gam)}\Big)^* \quad,
\end{alignat*}
This ODE has solutions
\begin{alignat*}{5}
    z &=  \rad e^{i\ta} \quad,\quad \rad(\ta) =  (2\pi\ta)^{-1} \quad,\quad \gam =  (2\pi\ta)^{1-2\mv} 
\end{alignat*}
where $\ta\in(0,\infty)$ is the angle traversed around the origin when starting at infinity, while $\rad$ is distance from
the spiral center. Such solutions parametrize algebraic spirals. For $\mv= 1$, the case most relevant for compressible flow,
hyperbolic spirals result. 

These approximations are considered correct at least near the spiral center. However, despite various attempts since Kaden's 1931 paper \cite{kaden-1931},
(see \cite{moore-1975} and references therein), existence of such solutions to the full Birkhoff-Rott equation has never been 
proven. Apart from applications, vortex spiral solutions also exhibit theoretically relevant phenomena such as apparent non-uniqueness
for the initial-value problem \cite{pullin-nuq,elling-hyp2010}, questions which are unlikely to be settled without 
rigorous existence proofs. 

\begin{figure}[h]
\includegraphics[width=.33\linewidth]{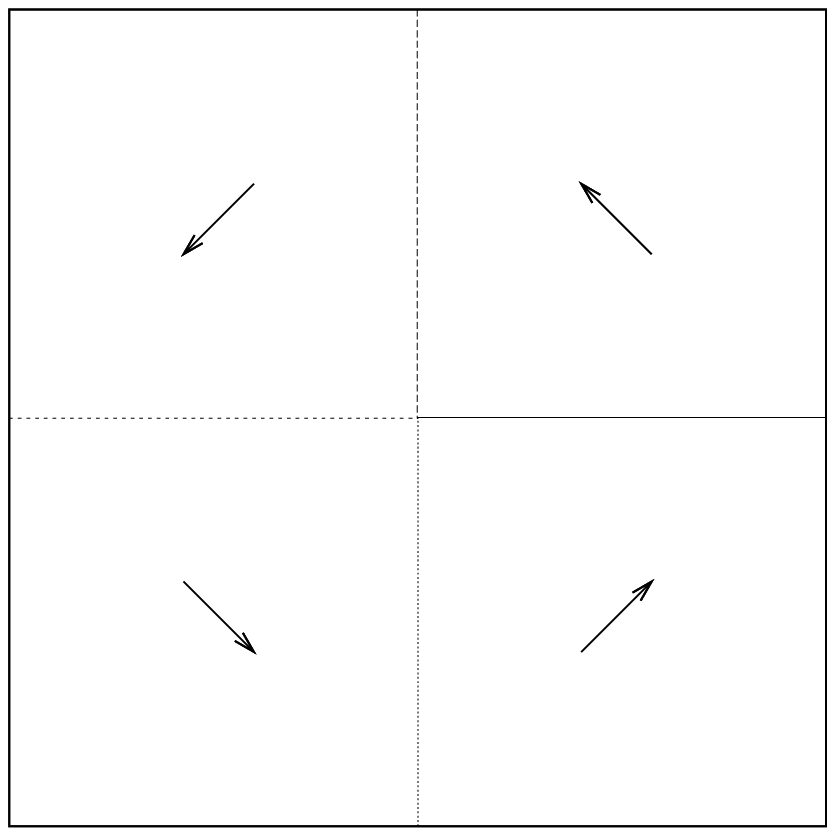}%
\includegraphics[width=.33\linewidth]{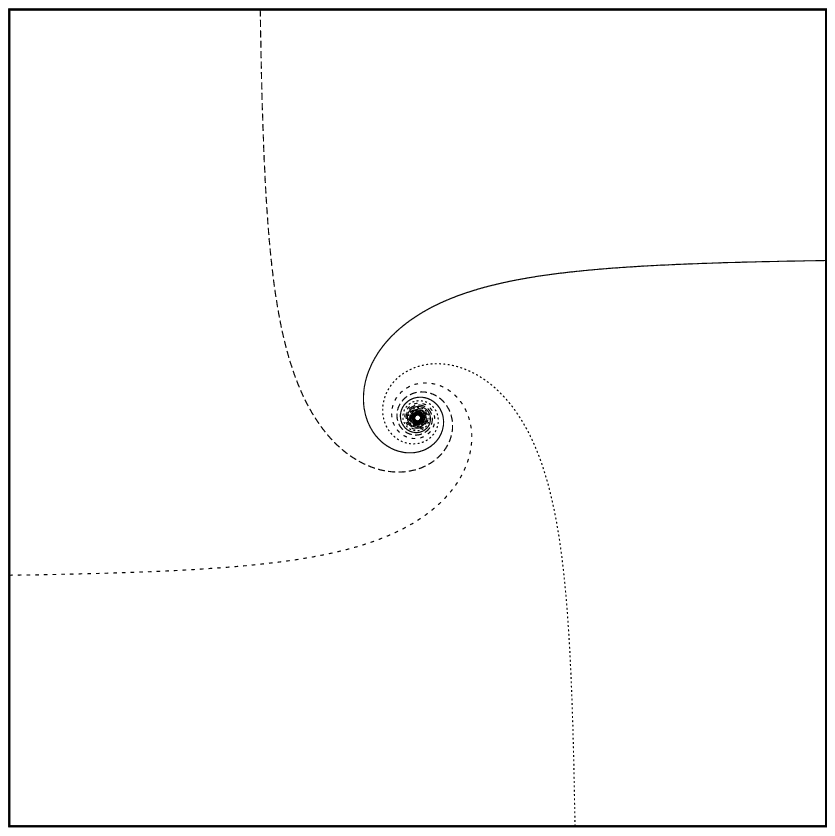}%
\includegraphics[width=.33\linewidth]{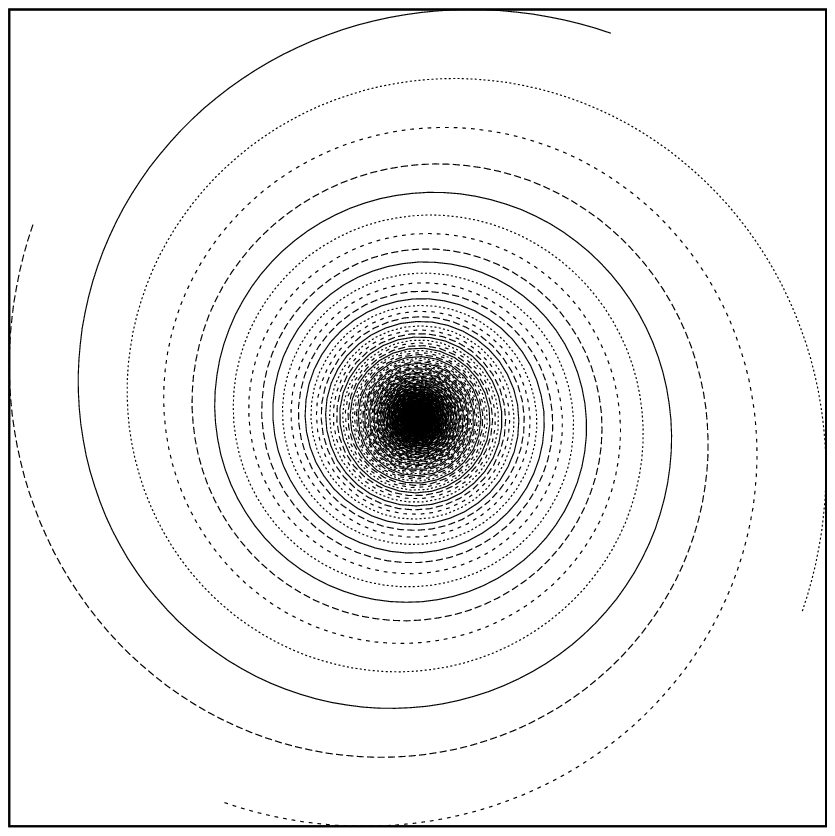}
\caption{Left: initial data ($\Nv= 4,\mv= 1$); center: rollup; right: core}
\mylabel{fig:Nspiral}
\end{figure}

In prior work \cite{elling-hyp2010,elling-Nspiral} the author constructed the first rigorous examples of algebraic vortex spiral solutions: 
in the case of $\Nv$
symmetric spiral branches with sufficiently large $\Nv$ (see \myref{fig:Nspiral}). 
However, working with the Birkhoff-Rott equation seems to lead to long and tedious proofs. Apart from 
working with principal value integral operators rather than 
derivatives, it appears that the $Z$ coordinates, 
while intuitively appealing, are not convenient for analysis. 

More importantly, the Birkhoff-Rott equation is relevant only to the extent that it yields particular solutions of the Euler
equations (see \cite{lopes-lopes-schochet} for a discussion of when this is true). 
It is natural to try to focus on the vorticity formulation \myeqref{eq:vort-div} of the Euler equations directly, as we do in this paper.
We prove existence of self-similar flows with \emph{smooth}\footnote{except in the origin} vorticity fields (rather than sheets) that stratify into algebraic spirals. 
These variants are also of independent interest, for example as models of eddies in turbulent flows. 

\subsection{Main result}

We consider solutions of the 2d incompressible Euler equations \myeqref{eq:euler} with the following properties: self-similarity
\begin{alignat*}{5}
    \vv(\vx,t) &= t^{\mv-1} \cvv(\cvx) \quad,\quad \vx =  t^\mv \cvx \quad,\quad \mv\in\boi{\frac12}\infty \myeqlab{eq:vx-cvx}
\end{alignat*}
(The exponent $\mv-1$ is obtained by dimensional analysis in section \myref{section:selfsimilarity}). 

The vorticity equation \myeqref{eq:vort-div} reduces to 
\begin{alignat*}{5} 0 &=  (\ssac\vv-\mv\cvx)\cdot\nabla_\cvx\cvort - \cvort \myeqlab{eq:ssvort} \end{alignat*}
where $\vort(\vx,t)=t^{-1}\cvort(\cvx)$. 
In analogy to steady flow, where integral curves of $\vv$ (appearing in the steady vorticity equation $\vv\cdot\nabla\vort=0$) are called \defm{streamlines},
we call $\ssac\vv-\mv\cvx$ \defm{pseudo-velocity} and its integral curves \defm{pseudo-streamlines}. 
\myeqref{eq:ssvort} suggests $\cvort$ will have higher regularity along those curves, which is observed physically as stratification;
self-similar vortex sheets coincide with pseudo-streamlines. 

We seek self-similar solutions that satisfy the initial condition 
\begin{alignat*}{5} \vort(\vx,t) \rightarrow \rad^\vrdiniex \vrdini(\pola) \quad\text{as}\quad t\searrow 0 \myeqlab{eq:vort-ini-cond} \end{alignat*} 
in some sense, where $(\rad,\pola)$ 
are polar coordinates centered in $\vx=0$,
and $\vrdiniex=-\frac1\mv$. 
(We show in section \myref{section:selfsimilarity} that $\mv$ is determined by choosing $\vrdiniex$, i.e.\ by choosing the initial data.) 

Since we are in two space dimensions, $\vort(\vx,t=0)$ is locally integrable if and only if $\vrdiniex=-\frac1\mv>-2\eqv\mv\in\boi{\frac12}{\infty}$, 
so we consider only this range of $\mv$.
The most relevant case, especially from the perspective of compressible flow, is $\mv=1$, the \defm{acoustic scaling},
which yields $\vv$ that are bounded in time and space. 
The famous Kaden spiral \cite{kaden-1931} is for $\mv=\frac23$. 

If initially there is net vorticity near the origin, i.e.\ if $\int_0^{2\pi}\vrdini(\pola)d\pola=\vrdini^\wedge(0)\neq 0$, 
then physical intuition suggests that for $t>0$ the flow rotates near the origin, so that the pseudo-streamlines
are distorted from initially straight rays into spirals. More precisely they will become \defm{algebraic} spirals: each streamline can be parametrized as
\begin{alignat*}{5} \pola \mapsto \rad(\pola) &= \pola^{-\mv} (C + o(1)) \quad\text{as}\quad \pola\rightarrow\infty \end{alignat*} 
This behaviour corresponds to the closely packed spirals seen in Figure \myref{fig:Nspiral}.

Let $\Wi:=\ell^1(\Z)^\vee$ be the Wiener algebra (functions with absolutely summable Fourier coefficients).
\begin{theorem}
    \mylabel{th:mainresult}%
    Let $\mv\in(\frac12,\infty)$. Let $\vrdinic\neq 0$ be constant.
    Assume:
    \begin{enumerate}
    \item Sufficiently high periodicity: $\vrdini$ is $\frac{2\pi}{\Nv}$-periodic for sufficiently large $\Nv\in\N$ 
    \item Smallness: $\vrdini=\vrdinic+\vrdinil$ for $\vrdinil\in\Wi$ with $\norma{\vrdinil/\vrdinic}{\Wi}$ sufficiently small.
    \end{enumerate}
    Then there is an $\stf$ with the aforementioned properties which defines
    a classical solution of the Euler equations \myeqref{eq:euler} for $(\vx,t)\in(\R^2\backslash\{0\})\times\roi0\infty$.
    For $\mv>\frac23$ it is a weak solution at $\vx=0$ as well. 
\end{theorem}

\begin{remark}
    We expect the solutions to be weak solutions for $\mv\leq\frac23$ as well; this will be analyzed in future work.
\end{remark}

\begin{remark}
    For $\vrdinic=0$ 
    spiral rollup should not be expected for arbitrary small perturbations; generically 
    and both the physical behaviour and mathematical analysis of the flow will be entirely different.
\end{remark}

\begin{remark}
    $\Wi$ includes\footnote{\cite[Theorem I 6.2]{katznelson}} in particular all H\"older-continuous functions with H\"older exponent $s>\frac12$, hence all $\Cinf$ functions. 
\end{remark}

\begin{remark}
    The restriction to sufficiently large $\Nv$ will be relaxed in later work, but it cannot be expected to be fully removable, 
    as simple physical arguments show.
    Consider for example $\mv=1$ and a single vortex sheet on the positive horizontal axis with uniform vorticity distribution:
    $Z(\Gamma,0)=\Gamma$ (no symmetry, corresponding to $\Nv=1$). 
    This would cause the velocity integral in the Birkhoff-Rott equation \myeqref{eq:br} to diverge
    due to the contribution at $\Gamma'\rightarrow\infty$, corresponding to spatial infinity. 
    An equal and opposite sheet on the negative horizontal axis (which would correspond to $\Nv=2$ symmetry) is needed to produce an
    opposite singularity that cancels the first one.
\end{remark}

\begin{remark}
    Global existence (i.e.\ without imposing smallness for $\vrdinil$) is not expected to hold:
    for the closely related vortex sheet spirals, the numerical work of Pullin \cite{pullin-jfluid-1978,pullin-nuq} 
    shows complicated bifurcation phenomena with limit points, pitchfork bifurcations etc., hence non-existence and
    non-uniqueness (which is our main motivation).
\end{remark}

\begin{remark}
    The algebraic spirals are not just mathematical constructions, but physically noticeable due to ``stratification'':
    near the spiral center consider line segments $I$ in radial direction, 
    starting and ending in the same pseudo-streamline which passes through a single spiral turn in between
    (hence $I$ crosses all the other pseudo-streamlines once). Then $\osc_I\vort=\sup_I\vort-\inf_I\vort\geq\delta\sup_I|\vort|$ 
    for some $\delta>0$ independent of $\rad,\pola$ as $\rad\searrow 0$, \emph{unless} $\vrdini$ is constant.
    (This is shown in Section \myref{section:stratification}.)
\end{remark}

\subsection{Overview}

\mylabel{section:overview}

The solutions we construct are perturbations of a ``trivial'' \defm{background solution}. To obtain the latter, we observe that 
\emph{radial} initial data $\vrdini(\pola)=\vrdinic=\text{const.}$ in \myeqref{eq:vort-ini-cond} extends to a self-similar solution
\begin{alignat*}{5} 
    \cvort(\cvx) &= \crad^{-\frac1\mv}\vrdinic \csep \crad=|\cvx| \quad.
\end{alignat*}
of the problem.
Then
\begin{alignat*}{5} 
    \vort(\tv,\vx) &= \tv^{-1}\cvort(\cvx) \topref{eq:vx-cvx}{=} \rad^{-\frac1\mv} \vrdinic \quad.
\end{alignat*}
(This background solution happens to be steady, but the perturbations we construct for non-constant $\vrdini$ are not steady.)
From the pseudo-velocity $\cvv-\mv\cvx$ we obtain its integral curves, the \defm{pseudo-streamlines} (see solid curves in Figure \myref{fig:ta-uan}), 
in the expected form
\begin{alignat*}{5} \crad = (\ta-\ta_0)^{-\mv} \crad_0  \csep \ta=\measuredangle\cvx \end{alignat*} 
for some constant $\crad_0=\crad_0(C)$ and arbitrary constant rotation $\ta_0$. 
Curves with this radius-angle relationship (in polar coordinates) are called \defm{algebraic spirals}; 
for $\mv=1$ the \defm{hyperbolic spirals} $\crad=\crad_0\frac{1}{\ta}$ are obtained. Such spirals feature densely packed turns in the inner
part, as seen in Figure \myref{fig:Nspiral} right. These turns will self-intersect when we make perturbations that do not decay sufficiently
rapidly (see Figure \myref{fig:self-intersection} right). Obtaining this decay is our biggest\footnote{Away from the origin and from spatial infinity, 
the problem consists of a uniformly elliptic equation coupled to a non-degenerate transport equation,
for which small perturbations would be relatively easy to obtain.} technical obstacle. 

In section \myref{section:prelim} we establish some notation and well-known results for easy reference. 
Section \myref{section:equations} performs a succession of coordinate changes, which are the key to the solution. 
The problem becomes tractable when expressed in coordinates $\Pav=(\ta,\uan)$,
where $\ta\in\boi0\infty$ is angle along a particular 
pseudo-streamline ($\ta=0+$ corresponds to spatial infinity; each increase of $\ta$ by $2\pi$ passes through another turn of the spiral and $\ta\rightarrow\infty$ approaches the spiral center), while
$\uan$ is the polar angle at $\crad\rightarrow\infty$ ($\uan$ parametrizes the family of pseudo-streamlines). 
This is a rather complicated implicitly defined change of coordinates which depends not only on $\cvx$ but also on the a priori unknown
solution $\stf$ (since we require $\pta$ to be tangential to its pseudo-streamlines). 
The complications are rewarded by some ``miraculous'' formulas that solve all equations except
the curl constraint \myeqref{eq:vort-def} which becomes a 3rd order nonlinear PDE. 

In section \myeqref{section:linearization} we linearize this PDE around the background solution above. After a Fourier transform
the PDE decouples into infinitely many 3rd order ODE in the Fourier variable $\etav\in\R$ dual to $\ta$.
The family of ODE is parametrized by the Fourier variable $\nv\in\Z$ dual to $\uan$. 
The ordinary differential operator $\Lopn=\Rop+\Eop$ has a dominant part $\Rop$ that factors into three Fuchsian 1st order operators
(see \myeqref{eq:linearization}). Having established their basic properties in section \myref{section:intops}, 
we show in section \myref{section:inversion} that $\Lopn$ is invertible if $\nv$ is zero or sufficiently large.
In section \myref{section:existence} we choose appropriate function spaces for $\stf$ and for the defect of the PDE, 
then apply the implicit function theorem to obtain solutions of the problem. The following sections analyze additional properties (regularity,
proper weak solutions, pseudo-streamline asymptotics, stratification). 

\section{Preliminaries}

\mylabel{section:prelim}

Here we establish some basic notation and elementary results.

\subsection{Notation}

For vectors $\va,\vb$, $\va\cdot\vb\defeq\va^T\vb$, $|\va|^2\defeq\va\cdot\va$ while $\va^2= \va\va^T$ is a rank 1 matrix. 
We use the Frobenius product $A:B:=\trace(A^TB)$.

All derivatives are meant in the distributional sense, whenever incompatible with the pointwise sense.
$\nabla$ is the gradient as a column vector, $\nabla^T$ as a row vector.
$\nabla^2=\nabla\nabla^T$ is the Hessian while $\Delta= |\nabla|^2=\nabla^T\nabla$ is the Laplacian.
\begin{alignat*}{5}
    \skw &\defeq \begin{bmatrix}
        0 & -1 \\
        1 & 0
    \end{bmatrix}
    \quad\impl\quad
    \skw^T= -\skw
    \quad,\quad
    \skw\skw= -I 
    \csep
    \skw^T\skw = \skw\skw^T = I 
    \myeqlab{eq:skw}
    \\
    \va^\perp &\defeq \skw\va
    =  \begin{bmatrix}
        0 & -1 \\
        1 & 0
    \end{bmatrix}
    \begin{bmatrix}
        a^1 \\
        a^2
    \end{bmatrix} 
    =  
    \begin{bmatrix}
        -a^2 \\
        a^1
    \end{bmatrix}
    \myeqlab{eq:perp-skw}
    \\
    \va \times \vb &\defeq \va^\perp\cdot\vb =  (\skw\va)^T\vb =  \va^T\skw^T\vb =  \begin{bmatrix}
        -a^2 & a^1 
    \end{bmatrix}\begin{bmatrix}
        b^1 \\
        b^2
    \end{bmatrix} 
    =  a^1b^2-a^2b^1
    =  \Det\begin{bmatrix}
        a^1 & b^1 \\
        a^2 & b^2
    \end{bmatrix} =  \Det\begin{bmatrix}
        \va & \vb
    \end{bmatrix}
    \myeqlab{eq:cross-skw}
\end{alignat*}
\begin{alignat*}{5}
    \nabla\times\vv &\defeq \nabla^\perp\cdot\vv
    =  \begin{bmatrix}
        -\partial_2 & \partial_1
    \end{bmatrix}\begin{bmatrix}
        v^1 \\
        v^2
    \end{bmatrix} 
    =  \partial_1v_2-\partial_2 v_1
\end{alignat*}
\begin{alignat*}{5}
    -\va^T\skw\va &=  \va^T\skw^T\va =  (\skw\va)^T\va =  \va\times\va =  a^1a^2-a^2a^1 =  0
    \\
    \Det\begin{bmatrix}
        \va & \va^\perp 
    \end{bmatrix}
    &= 
    \va^T\skw^T(\skw\va)
    = 
    |\va|^2
    \myeqlab{eq:det-vec-perp}
\end{alignat*}
$\adjt A$ is the algebraic adjoint (the transpose of the cofactor matrix).
\begin{alignat*}{5} 
    \begin{bmatrix}
        a_{22} & -a_{12} \\
        -a_{21} & a_{11}
    \end{bmatrix}
    \begin{bmatrix}
        0 & -1 \\
        1 & 0
    \end{bmatrix}
    &=
    \begin{bmatrix}
        -a_{12} & -a_{22} \\
        a_{11} & a_{21} 
    \end{bmatrix}
    =
    \begin{bmatrix}
        0 & -1 \\
        1 & 0
    \end{bmatrix}
    \begin{bmatrix}
        a_{11} & a_{21} \\
        a_{12} & a_{22}
    \end{bmatrix}
    \\\quad\impl\quad
    (\adjt A)\skw 
    &=
    \skw A^T
    \myeqlab{eq:adjt-skw}
\end{alignat*}

$\T=\roi0{2\pi}$ comes with the circle metric. Functions defined on $\T$ are meant to be $2\pi$-periodic (in that variable).

$f^\wedge(\etav):=(2\pi)^{-1}\int_\R f(x)e^{-i\etav x}dx$ is the Fourier transform of $f$ (in all variables),
$f^\vee(x):=\int_\R f(\etav)e^{i\etav x}d\etav$ the inverse Fourier transform, 
$f*g(\etav):=\int_\R f(\xiv)g(\etav-\xiv)d\xiv$ convolution. 
Then $(\partial_x f)^\wedge=i\etav f^\wedge$, $(x f)^\wedge=i\peta(f^\wedge)$, $1^\wedge=\delta$, $\delta*f=f$. 

The Fourier transform for $2\pi$-periodic $f$ is $f^\wedge(\nv):=(2\pi)^{-1}\int_0^{2\pi} f(\uan)e^{-i\nv\uan}d\uan$, 
with inverse $f^\vee(\uan)=\sum_{\nv\in\Z}f(\nv)e^{i\nv\uan}$
and convolution $f*g(\uan)=\int_{-\pi}^\pi f(\uan')g(\uan-\uan')d\uan'$. 

We write $\moac f$ or $\moacs f$ for the operator of pointwise multiplication by a function $f$, 
so that $(\partial_x\moacs f)g=\partial_x(fg)$ while $(\partial_xf)g=g\partial_xf$. 

We write $a\lesssim b$ if $\frac{a}{b}$ is defined and bounded as some limit is taken.
$a\sim b$ means $b\lesssim a$ as well.

$\asy x:=(1+|x|^2)^{\frac12}$.

\subsection{Banach spaces}

\subsubsection{Notation}

For an element $v$ of some vector space $V$ over a field $\K$, $\K v$ is the span of $v$.
$[X,Y]$ is the Banach space of linear continuous maps on $X$ into $Y$ with the uniform topology (operator norm); $[X]:=[X,Y]$.
$X\embed Y$ means $X$ is a subset of (or embedded in a canonical way in) $Y$, with stronger norm (not necessarily strictly stronger).

$\Ck{k}(X;Y)$ is the space of maps from $X$ to $Y$ that are $k$ times continuously differentiable.
$\Cbk{k}(X;Y)$ is the subspace of those maps whose derivatives of order up to $k$ are also bounded;
it is a Banach space with the usual norm. 
$\Cubk{k}(X;Y)$ is the closed subspace of functions whose derivatives of order up to $k$ are also \emph{uniformly} continuous.
$\Cuvk0(X;Y)$ is the closed subspace of $\Cubk0(X;Y)$ of functions that decay uniformly at infinity.

$\lone(\Z)$ is the space of absolutely summable sequences;
$\loner(\R)=\{(a_\nv)\in\lone(\Z)~|~\forall \nv\in\Z\setdiff(\Nv\Z):a_\nv=0\}$ is a closed subspace. (It will contain Fourier coefficients of functions that are not only $2\pi$-periodic but $\frac{2\pi}{\Nv}$-periodic.)
$\Ms$ is the space of measures of finite total variation on $\R$; $\Ms_\loc(\Omega)$ for $\Omega\subset\R$ open are distributions that are in $\Ms$ when restricted to compact subsets of $\Omega$. 
$\BV$ are functions of bounded variation (in the Tonelli-Cesari sense), 
$\BV_\loc(\Omega)$ are distributions that are $\BV$ when restricted to compact subsets of $\Omega$.
We write $\R$ (or $\CC$) for subspaces of constant functions.

\subsubsection{Induced norm}

Given a linear map $f:A\rightarrow B$ for vector spaces $A,B$ and given a Banach space $Y$ that is a linear subspace of $B$, 
$f^{-1}[Y]$ is a linear subspace of $A$ which carries the \defm{induced seminorm}  
\begin{alignat*}{5} \norma{x}{f^{-1}[Y]} \defeq \norma{f(x)}{Y}. \end{alignat*}
It is a normed space if $f$ is injective; it is Banach if $f(A)\subset Y$ is closed. 
Given a bijective map $g:X\rightarrow Y$ we define the induced norm on $g(X)$ using $g^{-1}:Y\rightarrow X$.

\subsubsection{Intersections and sums}

We define the canonical seminorm on $\bigcap_{i= 1}^m X_i$ as
\begin{alignat*}{5} \norma{a}{\bigcap_{i= 1}^m X_i} \defeq \sup_{i= 1}^m \norma{a}{X_i} \end{alignat*}
and the canonical seminorm on $\sum_{i= 1}^m X_i$ is
\begin{alignat*}{5} \norma{a}{\sum_{i= 1}^m X_i} \defeq \inf\Big\{ \sum_{i= 1}^m \norma{a_i}{X_i} ~\Big|~ a= \sum_{i= 1}^m a_i\ ,\ \forall i\in I:a_i\in X_i \Big\}\quad. \end{alignat*}
All Banach spaces we consider are continuously embedded in $\Ddist(\R\times\T)$, the space of distributions on $\R\times\T$. The topology of $\Ddist(\R\times\T)$ is Hausdorff, and the topology it induces on each $X_i$ is weaker than the $X_i$ norm topology. Therefore the seminorms above are norms and our intersections and sums are again Banach spaces.

\subsection{General change of coordinates}

When dealing with changes of coordinates we may 
use sub- or superscripts in round brackets to denote the representation in a particular frame. 
Transformation of divergences: say 
\begin{alignat*}{5} \nabla_\vx^T\vvxx = f\quad. \myeqlab{eq:divxx} \end{alignat*} 
Then for a compactly supported test function $\phi$
\begin{alignat*}{5} \int\phi f~d\vx
    = \int\phi f |\Det\nabt_\vy\vx| d\vy
\end{alignat*} 
while
\begin{alignat*}{5}&
    \int\phi f d\vx 
    = \int \phi \nabt_\vx\vvxx d\vx 
    =  -\int\nabt_\vx\phi\vvxx d\vx 
    =  -\int\nabt_\vy\phi\nabt_\vx\vy~\vvxx d\vx 
    \\&=  -\int\nabt_\vy\phi\nabt_\vx\vy~\vvxx |\Det\nabt_\vy\vx| d\vy
    =  \int \phi \nabla_\vy^T(\nabt_\vx\vy~\vvxx |\Det\nabt_\vy\vx| ) d\vy
\end{alignat*}
(the integrals are interpreted as duality for distributions). 
The last two right-hand sides are equal, and we may drop $|\cdot|$ since the signs cancel. By varying $\phi$ we obtain
\begin{alignat*}{5} 
    f = \nabt_\vx\vvxx \quad\eqv\quad
    f \Det\nabt_\vy\vx 
    &= \nabla_\vy^T (\Det\nabt_\vy\vx \nabt_\vx\vy~\vvxx ) 
    = \nabla_\vy^T (\Det\nabt_\vy\vx~\vv^{(\vy)}) 
    \\&= \nabla_\vy^T (\adjt\nabt_\vy\vx~\vvxx) 
    = \nabla_\vy^T (\frac{\nabt_\vx\vy}{\Det\nabt_\vx\vy}\vvxx) 
    \myeqlab{eq:div-xy}
\end{alignat*} 
Laplacian: $\Delta_\vx u=\nabt_\vx\nabla_\vx u= f$ transforms to 
\begin{alignat*}{5}
    f \Det\nabt_\vy\vx 
    &\topref{eq:div-xy}{=} 
    \nabt_\vy (\adjt \nabt_\vy\vx \nabla_\vx u) 
    = 
    \nabt_\vy (\Det\nabt_\vy\vx \nabt_\vx\vy (\nabt_\vx\vy)^T \nabla_\vy u) 
    \myeqlab{eq:Lap-xy}
\end{alignat*}

\subsection{Conformal change of coordinates}

Let $A\in\R^{n\times n}$. We say $A$ is \defm{conformal} if 
\begin{alignat*}{5} A^TA &=  sI \myeqlab{eq:ATA-sI} \end{alignat*} 
for some \defm{conformal factor} $s>0$. 
For $2\times 2$ conformal matrices $s^2=\Det(A^TA)=(\Det A)^2$, so $\Det A= \pm s$. 
An invertible matrix preserves orientation if and only if $\Det A>0$, so that 
\begin{alignat*}{5} \Det A &=  s . \myeqlab{eq:det-conformal} \end{alignat*} 
For a \emph{conformal orientation-preserving} change of coordinates (i.e.\ all Jacobians conformal and orientation-preserving)
\begin{alignat*}{5}
    \Delta_\vx\phi = f
    \quad\eqv\quad
    f \Det\nabt_\vy\vx
    &\topref{eq:Lap-xy}{= } 
    \nabt_\vy (
    \Det\nabla_\vy^T\vx
    \subeq{
        \nabt_\vx\vy(\nabt_\vx\vy)^T
    }{
        \toprefb{eq:ATA-sI}{eq:det-conformal}{= }
        \frac1{\Det\nabla_\vy^T\vx}I 
    }
    \nabla_\vy\phi)
    =  
    \Delta_\vy\phi 
    \myeqlab{eq:conformal-Delta}
\end{alignat*}

\section{Equations}

\mylabel{section:equations}

In this section we perform a succession of coordinate changes which are essential to making the problem solvable.

\subsection{Self-similarity}

\mylabel{section:selfsimilarity}

Here we perform a careful but elementary derivation of the self-similar vorticity equation.

We change from variables $\vy=(t,\vx)$ to $\cvy=(\ct,\cvx)$ with $\ct=t,\cvx\topref{eq:vx-cvx}{=}t^{-\mv}\vx$. Then
\begin{alignat*}{5} \nabt_\vy\cvy = \nabt_{(t,\vx)}\begin{bmatrix} \ct \\ \cvx \end{bmatrix}
    =
    \begin{bmatrix}
        1 & 0 \\
        -\mv t^{-\mv-1}\vx & t^{-\mv}I
    \end{bmatrix} 
    = 
    \begin{bmatrix}
        1 & 0 \\
        -\mv \ct^{-1}\cvx & \ct^{-\mv}I
    \end{bmatrix} 
    \myeqlab{eq:cvy-cvx} 
\end{alignat*} 
\begin{alignat*}{5} \Det\nabt_\vy\cvy \topref{eq:cvy-cvx}{=} \ct^{-\ndim\mv} 
    \myeqlab{eq:cvy-cvx-det} \end{alignat*}
Vorticity equation: 
\begin{alignat*}{5} 0 
    &\topref{eq:vort-div}{=}
    \nabt_\vy\begin{bmatrix}
        \vort \\
        \vort\vv
    \end{bmatrix}
    \topref{eq:div-xy}{=} 
    \nabt_\cvy \Big( \frac{\nabt_\vy\cvy}{\Det\nabt_\vy\cvy}\begin{bmatrix}
        \vort \\
        \vort\vv
    \end{bmatrix} \Big)
    \toprefb{eq:cvy-cvx}{eq:cvy-cvx-det}{=}
    \nabt_\cvy(
    \ct^{\ndim\mv}
    \begin{bmatrix}
        1 & 0 \\
        -\mv \ct^{-1}\cvx & \ct^{-\mv}I
    \end{bmatrix} 
    \begin{bmatrix}
        \vort \\
        \vort\vv
    \end{bmatrix} )
    \\&=
    \pct(\ct^{\ndim\mv}\vort) + \nabt_\cvx[\ct^{\ndim\mv}\vort(\ct^{-\mv}\vv-\mv \ct^{-1}\cvx)]
    \myeqlab{eq:cvort-pre}
\end{alignat*}
We seek \defm{self-similar} solutions
\begin{alignat*}{5} \vv(\tv,\vx)=\ct^{\mv-1}\cvv(\cvx) \myeqlab{eq:vv-cvv} \end{alignat*}
($\ct$ exponent $\mv-1$ because we cannot expect a large variety of such solutions unless $\vv,\cvx$ on the right-hand side of \myeqref{eq:cvort-pre} come with the same $\ct$ exponent). Then
\begin{alignat*}{5} 
    0 &\topref{eq:cvort-pre}{=} \pct(\ct^{\ndim\mv}\vort) + \nabt_\cvx[\ct^{\ndim\mv-1}\vort(\subeq{\cvv-\mv\cvx}{=:\cvz})]
    \myeqlab{eq:cvz}
\end{alignat*}
$\cvz$ is the \defm{pseudo-velocity}. 
Now
\begin{alignat*}{5} 
    \vort \topref{eq:vort-def}{=} \nabla_\vx\times\vv 
    =
    \ct^{-\mv} \nabla_\cvx\times(\ct^{\mv-1}\cvv)
    =
    \ct^{-1} \subeq{\nabla_\cvx\times\cvv}{\eqdef\cvort}
    \myeqlab{eq:vort-cvort}
\end{alignat*} 
and thus 
\begin{alignat*}{5} 0 \toprefb{eq:cvz}{eq:vort-cvort}{=} 
    \pct(\ct^{\ndim\mv-1}\cvort) + \nabt_\cvx(\ct^{\ndim\mv-2}\cvort\cvz)
    =
    (\ndim\mv-1)\ct^{\ndim\mv-2}\cvort + \ct^{\ndim\mv-1}\subeq{\pct\cvort}{=0} + \nabt_\cvx(\ct^{\ndim\mv-2}\cvort\cvz)
\end{alignat*} 
\begin{alignat*}{5} 
    \impl\quad
    0 =
    (\ndim\mv-1)\cvort + \nabt_\cvx(\cvort\cvz) \myeqlab{eq:cvort}
\end{alignat*} 
These steps are compatible with distributional interpretation of the outermost derivatives (since $(fg)'=f'g+fg'$ is valid for $f\in\Ddist$ and $g\in\Cinf$,
and $\ct^\beta$ is $\Cinf$ for $\ct>0$ and $\beta\in\R$).
Finally, 
\begin{alignat*}{5} 
    \cvv(\cvx) 
    &\topref{eq:vv-cvv}{=} \ct^{1-\mv}\vv(t,\vx) 
    \topref{eq:stf}{=} \ct^{1-\mv}\nabla_\vx^\perp\stf(t,\vx)
    \topref{eq:perp-skw}{=} \ct^{1-\mv} \skw\nabla_\vx\stf(t,\vx) 
    = \ct^{1-\mv} \skw\subeq{\nabla_\vx\cvx^T}{\topref{eq:vx-cvx}{=}\ct^{-\mv}I}\nabla_\cvx\stf(t,\vx) 
    \\&= \skw\nabla_\cvx[\ct^{1-2\mv}\stf(t,\vx)]
\end{alignat*} 
so since the left-hand side is constant in $t$, we have for some functions $f,g$ that
\begin{alignat*}{5} 
    t^{1-2\mv} \stf(t,\vx) = f(\cvx) + g(t)
\end{alignat*} 
and we may omit the irrelevant $g(t)$ and take
\begin{alignat*}{5} t^{1-2\mv} \stf(t,\vx) \eqdef \cstf(\cvx) \myeqlab{eq:stf-cstf} \end{alignat*} 
and then
\begin{alignat*}{5} \cvv = \skw\nabla_\cvx\cstf . \myeqlab{eq:cvv-cstf} \end{alignat*}

To see why $\vrdiniex=-\frac1\mv$ is the correct $\rad$ exponent in the initial condition \myeqref{eq:vort-ini-cond},
we make the following plausibility argument: for $\tv>0$ small, 
\begin{alignat*}{5} 
    \cvort(\cvx) 
    \topref{eq:vort-cvort}{=} 
    \tv \vort(\vx,\tv) 
    \approx 
    \tv \vort(\vx,0) 
    \topref{eq:vort-ini-cond}{=} 
    \tv \rad^\vrdiniex \vrdini(\pola) 
    \topref{eq:vx-cvx}= 
    \tv\cdot \tv^{\mv\vrdiniex} \crad^\vrdiniex \vrdini(\pola) 
    \myeqlab{eq:vrdiniex-calc}
\end{alignat*} 
Since the left- and right-hand side are (note $\pola=\measuredangle\vx=\measuredangle\cvx$) independent of $\tv$ except for $\tv^{1+\mv\vrdiniex}$,
we generally need $\vrdiniex=-\frac1\mv$.

\subsection{Change to conformal polar coordinates $\lrv,\pola$}

To study spirals converging to a common origin it is convenient to use some form 
of polar coordinates $(\crad,\pola)$, i.e.
\begin{alignat*}{5} \crad=|\cvx| \csep \pola=\measuredangle\cvx  \myeqlab{eq:crad-pola-def} \end{alignat*} 
(see Figure \myref{fig:lrv-pola}). We make the conformal and
orientation-preserving change to 
coordinates $\Tav= (\lrv,\pola)$ where 
\begin{alignat*}{5} e^{\lrv} = \crad \myeqlab{eq:rad-lrv} \end{alignat*} 
\begin{alignat*}{5}
    \cvx &=  (\cx,\cy) =  (\crad\cos\pola,\crad\sin\pola) =  e^\lrv(\cos\pola,\sin\pola) \\
    \pT \cvx &=  
    e^{\lrv}(\cos\pola,\sin\pola)= 
    \cvx \\
    \pa \cvx &=  
    e^{\lrv}(-\sin\pola,\cos\pola)= 
    \cvx^\perp \\
    \nabla_\Tav^T\cvx 
    &=  
    \nabt_{(\lrv,\pola)}\cvx
    =
    \begin{bmatrix}
        \cvx & \cvx^\perp 
    \end{bmatrix} 
    =  
    \begin{bmatrix}
        \cx & -\cy \\
        \cy & \cx
    \end{bmatrix} 
    \myeqlab{eq:xT}
    \\
    \adjt\nabt_\Tav\cvx 
    &=
    \begin{bmatrix}
        \cx & \cy \\
        -\cy & \cx 
    \end{bmatrix}
    =
    \begin{bmatrix}
        \cvx & \cvx^\perp
    \end{bmatrix}^T
    \myeqlab{eq:adjt-xT}
    \\
    \Det\nabla_\Tav^T\cvx &\topref{eq:det-vec-perp}{=} |\cvx|^2 =  \crad^2 
    =  e^{2\lrv}
    \myeqlab{eq:xT-det}
    \\
    \nabla_\cvx^T\Tav 
    &=  
    (\nabla_\Tav^T\cvx)^{-1} 
    =
    (\Det\nabla_\Tav^T\cvx)^{-1} 
    \adjt\nabla_\Tav^T\cvx
    =  \crad^{-2} \begin{bmatrix}
        \cx & \cy \\
        -\cy & \cx
    \end{bmatrix}
    =
    \crad^{-2} \begin{bmatrix}
        \cvx & \cvx^\perp
    \end{bmatrix}^T
    \myeqlab{eq:DDxa}
\end{alignat*}
Curl equation:
\begin{alignat*}{5}
    \cvort 
    &=  
    \Delta_\cvx\cstf 
    \quad\topref{eq:conformal-Delta}{\eqv }\quad
    \Det\nabt_\Tav\cvx~\cvort 
    \topref{eq:xT-det}{= } 
    e^{2\lrv} \Delta_\Tav \cstf
    \myeqlab{eq:Delta-tav}
\end{alignat*}

\subsection{Change to pseudo-streamline angular coordinates $\ta,\uan$}

\subsubsection{Motivation and derivation}

\begin{figure}
\parbox[t]{1.9in}{
\input{lrv-pola.pstex_t}
\caption{$\Tav=(\lrv,\pola)$ coordinates}
\mylabel{fig:lrv-pola}}\parbox[t]{1.9in}{
\input{ta-uan.pstex_t}
\caption{$\Pav=(\ta,\uan)$ coordinates}
\mylabel{fig:ta-uan}}\parbox[t]{1.9in}{
\input{polb-Uan.pstex_t}
\caption{$(\polb,\Uan)$ coordinates}
\mylabel{fig:polb-Uan}}
\end{figure}
It is necessary to align one coordinate with the pseudo-streamlines, a technique commonly used for the 2d Euler equations or 
generally when transport equations are coupled with other PDE. 
To this end we change to coordinates $\Pav= (\ta,\uan)$ 
(see Figure \myref{fig:ta-uan})
which were already discussed in the context of the theorem. A key observation of our paper is that choosing an \emph{angle} $\ta$ 
as the along-pseudostreamline parameter 
leads to surprisingly simple expressions for $\lrv$ and $\cvort$.

We take \begin{alignat*}{5} \pola=\ta+\uan \myeqlab{eq:pola} \end{alignat*} 
and leave $\lrv=\lrv(\ta,\uan)$ general for now.
    
We use the convenient notation 
\begin{alignat*}{5} \pU \defeq \pu-\pta \myeqlab{eq:pU-def} \end{alignat*} 
($\pU$ would result from a further\footnote{We do not use these coordinates since $\polb,\Uan$ are naturally periodic, but then the change would not be bijective.} change of coordinates to $(\polb,\Uan)$ where
$\polb\defeq\ta+\uan= \pola$ and $\Uan\defeq\uan$; see Figure \myref{fig:polb-Uan}).
Then
\begin{alignat*}{5}
    \nabt_\Pav\Tav 
    =
    \nabt_{(\ta,\uan)}\begin{bmatrix}
        \lrv \\
        \pola
    \end{bmatrix}
    &= 
    \begin{bmatrix}
        \lrv_\ta & \lrv_\uan \\
        \pola_\ta & \pola_\uan
    \end{bmatrix} 
    \topref{eq:pola}{=}
    \begin{bmatrix}
        \lrv_\ta & \lrv_\uan \\
        1 & 1
    \end{bmatrix} 
    \myeqlab{eq:TP}
    \\
    \Det\nabla_\Pav^T\Tav &\topref{eq:TP}{=}  \lrv_\ta-\lrv_\uan \topref{eq:pU-def}{= } -\lrv_\Uan
    \myeqlab{eq:detTP}
    \\
    \adjt\nabla_\Pav^T\Tav 
    &\topref{eq:TP}{=}  \begin{bmatrix}
        1 & -\lrv_\uan \\
        -1 & \lrv_\ta
    \end{bmatrix} 
    \myeqlab{eq:TP-adjt}
    \\
    \begin{bmatrix}
        \ta_\lrv & \ta_\pola \\
        \uan_\lrv & \uan_\pola
    \end{bmatrix}
    =  \nabla_\Tav^T\Pav 
    =  (\nabla_\Pav^T\Tav)^{-1}
    = (\Det\nabt_\Pav\Tav)^{-1}\adjt\nabt_\Pav\Tav
    &\toprefb{eq:detTP}{eq:TP-adjt}{=} 
    \frac{1}{\lrv_\Uan} \begin{bmatrix}
        -1 & \lrv_\uan \\
        1 & -\lrv_\ta
    \end{bmatrix} 
    \myeqlab{eq:PT}
\end{alignat*}
\begin{alignat*}{5} 
    \pu\moac{\lrv_\ta}-\pta\moac{\lrv_\uan}
    &= \pu\moac{\lrv_\ta}-\pu\moac{\lrv_\uan}+\pu\moac{\lrv_\uan}-\pta\moac{\lrv_\uan}
    \\&\topref{eq:pU-def}{=} -\pu\moac{\lrv_\Uan} + \pU\moac{\lrv_\uan}
    \myeqlab{eq:puta-puU}
\end{alignat*} 
and analogously
\begin{alignat*}{5} \lrv_\ta\pu-\lrv_\uan\pta
    &= 
    \lrv_\uan\pU - \lrv_\Uan\pu
    \myeqlab{eq:puta-puU-nondiv}
\end{alignat*} 

\begin{alignat*}{5}
    \begin{bmatrix}
        \partial_\lrv \\ \partial_\pola
    \end{bmatrix}
    = 
    \nabla_{\Tav}
    &= 
    \nabla_\Tav\Pav^T
    \nabla_{\Pav}
    \topref{eq:PT}{=}
    \frac1{\lrv_\Uan}
    \begin{bmatrix}
        -1 & 1 \\
        \lrv_\uan & -\lrv_\ta
    \end{bmatrix} 
    \begin{bmatrix}
        \pta \\ \pu
    \end{bmatrix}
    = 
    \frac1{\lrv_\Uan}
    \begin{bmatrix}
        \pu-\pta \\ \lrv_\uan\pta-\lrv_\ta\pu
    \end{bmatrix}
    \toprefb{eq:pU-def}{eq:puta-puU-nondiv}{=}
    \frac1{\lrv_\Uan}
    \begin{bmatrix}
        \pU \\ \lrv_\Uan\pu-\lrv_\uan\pU 
    \end{bmatrix}
    \\&=
    \begin{bmatrix}
        \lrv_\Uan^{-1}\pU \\ \pu-\lrv_\uan\lrv_\Uan^{-1}\pU 
    \end{bmatrix}
    \myeqlab{eq:pTpa-ptapu}
\end{alignat*}

\subsubsection{$\pta$ along pseudo-streamlines}

Now we complete the definition of $\Pav$ by imposing that the pseudo-streamlines are tangential to the $\pta$ direction: $\cvz$, the tangent of pseudo-streamlines, when expressed as $\cvz^{(\Pav)}$, i.e.\ in the $\Pav=(\ta,\uan)$ frame, has $\uan$ component zero: first,

\begin{alignat*}{5} 
    \adjt\nabt_\Pav\Tav \adjt\nabt_\Tav\cvx~\cvv^{(\cvx)}
    \topref{eq:cvv-cstf}{=}
    \adjt\nabt_\Pav\Tav \adjt\nabt_\Tav\cvx \skw \nabla^{}_\cvx\cstf
    \topref{eq:adjt-skw}{=}
    \skw \nabla^{}_\Pav\Tav^T \nabla^{}_\Tav\cvx^T \nabla^{}_\cvx\cstf
    =
    \skw \nabla^{}_\Pav\cstf
    \myeqlab{eq:cstf-cvx-Pav}
\end{alignat*} 

\begin{alignat*}{5} 
    \adjt\nabt_\Pav\Tav~
    \frac{\adjt\nabt_\Tav\cvx}{\Det\nabt_\Tav\cvx}
    \cvx
    \toprefc{eq:TP-adjt}{eq:adjt-xT}{eq:xT-det}{=}
    \begin{bmatrix}
        1 & -\lrv_\uan \\
        -1 & \lrv_\ta
    \end{bmatrix}
    \crad^{-2}
    \begin{bmatrix}
        \cvx & \cvx^\perp
    \end{bmatrix}^T
    \cvx
    =
    \begin{bmatrix}
        1 & -\lrv_\uan \\
        -1 & \lrv_\ta
    \end{bmatrix}
    \begin{bmatrix}
        1 \\
        0 
    \end{bmatrix}
    =
    \begin{bmatrix}
        1 \\
        -1
    \end{bmatrix}
    \myeqlab{eq:cvx-part}
\end{alignat*} 

\begin{alignat*}{5} &
    \Det\nabt_\Tav\cvx
    \Det\nabt_\Pav\Tav
    ~\cvz^{(\Pav)}
    =
    \Det\nabt_\Tav\cvx
    \Det\nabt_\Pav\Tav
    \nabt_\Tav\Pav
    \nabt_\cvx\Tav
    ~\cvz^{(\cvx)}
    \topref{eq:cvz}{=}
    \adjt\nabt_\Pav\Tav~
    \adjt\nabt_\Tav\cvx
    (\cvv^{(\cvx)}-\mv\cvx)
    \\&\toprefb{eq:cstf-cvx-Pav}{eq:cvx-part}{=}
    \skw\nabla_\Pav\cstf
    -\mv
    \Det\nabt_\Tav\cvx
    \begin{bmatrix}
        1 \\
        -1 
    \end{bmatrix}
    \topref{eq:xT-det}{=}
    \begin{bmatrix}
        0 & -1 \\
        1 & 0
    \end{bmatrix}\begin{bmatrix}
        \cstf_\ta \\
        \cstf_\uan
    \end{bmatrix}
    - \mv 
    \crad^2
    \begin{bmatrix}
        1 \\
        -1 
    \end{bmatrix}
    \topref{eq:rad-lrv}{=}
    \begin{bmatrix}
        -\cstf_\uan \\
        \cstf_\ta 
    \end{bmatrix}
    + \mv 
    e^{2\lrv}
    \begin{bmatrix}
        -1 \\
        1 
    \end{bmatrix}
    \myeqlab{eq:cvz-Pav}
\end{alignat*} 
We want the $\uan$ (second) component of $\cvz^{(\Pav)}=\cvz^{(\ta,\uan)}$ zero:
\begin{alignat*}{5} 0 
    = 
    \cstf_\ta
    +\mv e^{2\lrv}
    \quad\eqv\quad 
    \crad^2 = e^{2\lrv} = \frac{\cstf_\ta}{-\mv} 
    \quad\eqv\quad  
    \text{\fbox{$\lrv = \frac12\log\frac{\cstf_\ta}{-\mv}$}}
    \myeqlab{eq:lrv-cstf}
\end{alignat*} 
Our solutions will be chosen so that $-\cstf_\ta/\mv>0$. Then
\begin{alignat*}{5} \lrv_x &\topref{eq:lrv-cstf}{= } \frac{\stf_{x\ta}}{2\stf_\ta}\quad(x= \ta,\uan,\Uan) \myeqlab{eq:lrv-x} \end{alignat*}
\begin{alignat*}{5} 
    \Det\nabla_\Pav^T\Tav 
    &\topref{eq:detTP}{=}
    -\lrv_\Uan
    \topref{eq:lrv-x}{=}
    -\frac{\stf_{\ta\Uan}}{2\stf_\ta}
    \myeqlab{eq:detTPb}
\end{alignat*} 
Returning to the expression for $\cvz$ we have
\begin{alignat*}{5} 
    \Det\nabt_\Tav\cvx
    \Det\nabt_\Pav\Tav
    ~\cvz^{(\Pav)}
    &\toprefb{eq:cvz-Pav}{eq:lrv-cstf}{=}
    \begin{bmatrix}
        -\cstf_\uan \\
        \cstf_\ta 
    \end{bmatrix}
    +
    \cstf_\ta
    \begin{bmatrix}
        1 \\
        -1 
    \end{bmatrix}
    \topref{eq:pU-def}{=}
    \begin{bmatrix}
        -\cstf_\Uan \\
        0
    \end{bmatrix}
    \myeqlab{eq:cvz-Pav-b}
\end{alignat*} 
Then the vorticity equation reduces to
\begin{alignat*}{5} 
    0 
    &\topref{eq:cvort}{=} 
    (\ndim\mv-1)\cvort + \nabt_\cvx(\cvort\cvz^{(\vx)}) 
    \topref{eq:div-xy}{=} 
    \Det\nabt_\Tav\cvx~\Det\nabt_\Pav\Tav (2\mv-1) \cvort + \nabt_\Pav [ \cvort\Det\nabt_\Tav\cvx\Det\nabt_\Pav\Tav ~\cvz^{(\Pav)} ]
    \\&\toprefc{eq:xT-det}{eq:detTP}{eq:cvz-Pav-b}{=}
    \crad^2 (-\lrv_\Uan) (2\mv-1) \cvort + 
    \begin{bmatrix}
        \pta \\
        \pu
    \end{bmatrix}^T 
    \Big( \cvort \begin{bmatrix}
        -\cstf_\Uan \\
        0 
    \end{bmatrix} \Big)
    \toprefb{eq:lrv-cstf}{eq:detTPb}{=}
    \frac{\cstf_\ta}{-\mv} 
    (- \frac{\cstf_{\Uan\ta}}{2\cstf_\ta})
    (2\mv-1) \cvort
    - \pta(\cvort\cstf_\Uan)
    \\&=
    (1-\frac1{2\mv}) \cstf_{\Uan\ta} \cvort - \pta(\cvort\cstf_\Uan)
\end{alignat*} 
This is solved by the surprisingly elegant formula
\begin{alignat*}{5} 
    \text{\fbox{$\cvort = \subeq{(\cstf_\Uan)^{-\frac1{2\mv}}}{=:\cvorf} \Vort(\uan) $}}
    \myeqlab{eq:vort-stf-ta}
\end{alignat*} 
where $\Vort$ is some function that can be chosen freely, as data.
That is natural: we have a choice of initial data that is also self-similar, hence constant along rays, which
are constant-$\uan$ contours. (Note that $\Vort$ does not correspond exactly to initial data; the link will be made in section
\myref{section:ini-data}.) 

It seems very fortunate to obtain the simple formulas \myeqref{eq:lrv-cstf} for $\lrv$ and \myeqref{eq:vort-stf-ta} for $\cvort$.
While the one for $\lrv$ is less surprising (it can be understood as pseudo-velocity $\cvz$ being orthogonal to pseudo-streamlines which are tangential to $\pta$), the one for $\cvort$ seems somewhat miraculous. ``Miracles'' in mathematics indicate a hidden underlying structure,
but it is not clear what it might be. In any case, it is not clear that this structure would survive generalization,
such as to compressible flow (where stream function formulations are much more awkward anyway), so we do not explore this point further.

\subsubsection{Curl equation (divergence form)}

The last remaining equation is the curl equation \myeqref{eq:Delta-tav}:
\begin{alignat*}{5}
    0 
    \topref{eq:Delta-tav}{=} 
    \Delta_\Tav\cstf-e^{2\lrv}\cvort 
    \quad\topref{eq:Lap-xy}{\eqv}\quad
    0 
    =
    \nabt_\Pav (\adjt\nabt_\Pav\Tav~\nabla_\Tav\cstf) - \Det\nabt_\Pav\Tav~e^{2\lrv} \cvort 
\end{alignat*}
Then use 
\begin{alignat*}{5}
    \nabla_\Pav^T(\adjt\nabt_\Pav\Tav\nabla_\Tav\cstf)
    &\topref{eq:TP-adjt}{= }
    \begin{bmatrix}
        \pta \\ \pu
    \end{bmatrix}^T(
    \begin{bmatrix}
        1 & -\lrv_\uan \\
        -1 & \lrv_\ta
    \end{bmatrix} 
    \begin{bmatrix}
        \pT \\
        \pa
    \end{bmatrix}
    \cstf)
    \topref{eq:pU-def}= 
    \begin{bmatrix}
        -\pU \\
        \pu\moac{\lrv_\ta}-\pta\moac{\lrv_\uan}
    \end{bmatrix}^T
    \begin{bmatrix}
        \pT \\
        \pa
    \end{bmatrix}
    \cstf
    \\&\topref{eq:puta-puU}= 
    \begin{bmatrix}
        -\pU \\
        \pU\moac{\lrv_\uan}-\pu\moac{\lrv_\Uan}
    \end{bmatrix}^T
    \begin{bmatrix}
        \pT \\
        \pa
    \end{bmatrix}
    \cstf
    =
    \big(-\pU[\pT-\moac{\lrv_\uan}\pa]-\pu\moac{\lrv_\Uan}\pa\big)\cstf
\end{alignat*}
and
\begin{alignat*}{5} 
    -\Det\nabt_\Pav\Tav~e^{2\lrv}\cvort 
    \toprefb{eq:detTPb}{eq:lrv-cstf}{=} 
    \frac{\cstf_{\ta\Uan}}{2\cstf_\ta} \frac{\cstf_\ta}{-\mv} \cvort
    =  
    -\frac{1}{2\mv} \cstf_{\ta\Uan} \cvort
\end{alignat*} 
to obtain
\begin{alignat*}{5}
    0
    =  
    - \pU\subeq{(\pT-\moac{\lrv_\uan}\pa)\cstf}{\eqdef\gU}
    - \pu \subeq{\moac{\lrv_\Uan}\pa\cstf}{\eqdef\gu}
    \subeq{-\frac1{2\mv} \cstf_{\ta\Uan} \cvort}{\eqdef\vortterm}
    \myeqlab{eq:Delta-pav-new}
\end{alignat*} 

\section{Linearization}

\mylabel{section:linearization}
\mylabel{section:scaling}\mylabel{section:background-values}

It is known from the numerical calculations
of Pullin \cite{pullin-jfluid-1978,pullin-nuq} that the physically closely related self-similar Birkhoff-Rott equation
has complicated bifurcation phenomena. Hence we cannot expect global existence and uniqueness for the more general 
self-similar Euler equations either. However, it is possible to obtain small-data results, by linearizing around the trivial background
solution (choose $\linv=0$ below) whose properties were outlined in section \myref{section:overview}. 

For linearization it is convenient to scale the background solutions from $\cstf\sim\ta^{1-2\mv}$ to $\sstf\sim 1$. 
In fact choosing exactly which power of $\ta$ to attach to a particular derivative is essential for solving our problem. 
In the final result, \myeqref{eq:linearization}, there will be one $\ta$ attached to $\pta,\pU$, while $\pu$ comes without $\ta$
(there is another choice since $\pU=\pu-\pta$). 
The careful distribution of $\ta$ over the terms and derivatives has an interesting effect: the asymptotic limits of spatial infinity 
($\crad=\infty$, corresponding to $\ta\searrow 0$) and of spiral center ($\crad\searrow 0$, corresponding to $\ta\uparrow\infty$) 
are treated simultaneously, even though their properties appear rather different in Figure \myref{fig:Nspiral}: 
the solutions have pseudo-streamlines that are straight lines at infinity, algebraic spirals near the center. 
The problem has a hidden homogeneity.

We define
\begin{alignat*}{5}
    \cstf &=  \ta^{1-2\mv} \sstf \quad,
    \myeqlab{eq:cstf-sstf}
    \\
    \sstf &=  \sstfz + \linv~\lsstf  \quad,
    \\
    \sstfz &=  \frac1{2\mv-1}\quad,
    \myeqlab{eq:sstfz}
    \\
    \sVort &= \Vortz + \linv~\Vortl   \quad,
    \myeqlab{eq:Vort-sVort}
    \\
    \Vortz &=  \frac{2\mv-1}{\mv}\quad.
    \myeqlab{eq:Vortz}
\end{alignat*}
where $\linv= 0$ will yield the background solution. 
For arbitrary $f$ abbreviate 
$\zeroacc f:=f_{|\linv=0}$ (values of an expression $f$ at the background solution)
and $\liac f:=(\partial_\linv f)_{|\linv=0}$ (first variation of $f$). 
We may regard $\linv$ as an independent variable like $\ta,\uan$ and commute $\partial_\linv$ with $\pu,\pta$ etc. 

The following scaling and linearization steps are tedious and elementary, so the reader may wish to skip 
them in the first pass.

\subsection{Scaled derivatives}

\begin{alignat*}{5}
    \mdta &\defeq \tao^{2\mv-1}\tao\pta\tao^{1-2\mv} =  \tao\pta+\moacs1-\moacs2\moacs\mv \myeqlab{eq:mdta} 
    \\
    \mdU 
    &\defeq 
    \tao^{2\mv-1}\tao\pU\tao^{1-2\mv} 
    \topref{eq:pU-def}{=} 
    \tao^{2\mv-1}\tao(\pu-\pta)\tao^{1-2\mv} 
    =  
    \tao\pu-(\tao\pta+\moacs1-\moacs2\moacs\mv)
    =
    \tao\pU+\moacs2\moacs\mv-\moacs1
    \myeqlab{eq:mdU} 
\end{alignat*}
We use without further mention that 
\begin{alignat*}{5} \zeroaccl{\modac\partial_{x_1}...\modac\partial_{x_m}\sstf}=\modac\partial_{x_1}...\modac\partial_{x_m}\sstfz 
    \myeqlab{eq:L-zeroacc}
\end{alignat*}
and
\begin{alignat*}{5} \liac{\modac\partial_{x_1}...\modac\partial_{x_m}\sstf}=\modac\partial_{x_1}...\modac\partial_{x_m}\lsstf 
    \myeqlab{eq:L-liac}
\end{alignat*}
for $x_1,...,x_m\in\{\ta,\uan,\Uan\}$. (This is not true for $\mda,\mdT$ etc., defined later, which are nonlinear in $\cstf$.)
Also, 
\begin{alignat*}{5} \zeroaccl{f\cdot g}=\zeroacc f\cdot\zeroacc g \quad, \end{alignat*} 
the product rule in the form
\begin{alignat*}{5} \liac{f\cdot g}=\zeroacc f\cdot\liac g+\liac f\cdot\zeroacc g \end{alignat*} 
and all other calculus rules apply.

\subsection{Linear operators}

\begin{alignat*}{5}
    \cstf_\uan 
    &\topref{eq:cstf-sstf}{=}
    \pu\tao^{1-2\mv}\sstf 
    =
    \tao^{1-2\mv}\mdu\sstf 
    \myeqlab{eq:cstf-u-sstf}
    \\
    \mdu\sstfz 
    &\topref{eq:sstfz}{=}
    \pu\frac{1}{2\mv-1} 
    =  
    0 
    \myeqlab{eq:mdu-sstfz} 
\end{alignat*} 

\begin{alignat*}{5}
    \cstf_\ta 
    &\topref{eq:cstf-sstf}{=}  
    \tao^{-2\mv}\tao^{2\mv-1}\tao\pta\tao^{1-2\mv}\sstf 
    \topref{eq:mdta}{=}  
    \tao^{-2\mv}\mdta\sstf 
    \myeqlab{eq:cstf-ta-sstf}
    \\
    \mdta\sstfz 
    &\toprefb{eq:mdta}{eq:sstfz}{=}  
    (\tao\pta+\moacs1-\moacs2\moacs\mv)\frac{1}{2\mv-1} 
    =  
    -1 
    \myeqlab{eq:mdta-sstfz} 
\end{alignat*} 

\begin{alignat*}{5}
    \cstf_\Uan 
    &\topref{eq:cstf-sstf}{=}  
    \tao^{-2\mv}\tao^{2\mv-1}\tao\pU\tao^{1-2\mv}\sstf 
    \topref{eq:mdU}{=}
    \tao^{-2\mv}\mdU\sstf 
    \myeqlab{eq:cstf-U-sstf}
    \\
    \mdU\sstfz 
    &\toprefb{eq:mdU}{eq:sstfz}{=}  
    (\tao\pU+\moacs2\moacs\mv-\moacs1)\frac{1}{2\mv-1} 
    =  
    1 
    \myeqlab{eq:mdU-sstfz} 
\end{alignat*} 

\begin{alignat*}{5}
    \cstf_{\uan\ta}
    &\topref{eq:cstf-sstf}{=}
    \tao^{-2\mv}\pu\tao^{2\mv-1}\tao\pta\tao^{1-2\mv}\sstf 
    \topref{eq:mdta}{=}
    \tao^{-2\mv}\mdu\mdta\sstf 
    \myeqlab{eq:cstf-ta-uan-sstf}
    \\
    \mdu\mdta\sstfz 
    &\topref{eq:mdta-sstfz}{=} 
    \pu(-1) 
    =  
    0  
    \myeqlab{eq:mdu-mdta-sstfz} 
\end{alignat*} 

\begin{alignat*}{5} \tao^2\pU\pta 
    &\topref{eq:pU-def}{=}
    \tao\tao(\pu-\pta)\pta
    =
    \tao[\pu\tao-(\pta\tao-\moacs 1)]\pta
    \topref{eq:pU-def}{=}
    (\tao\pU+\moacs 1)\tao\pta
    \myeqlab{eq:pU-pta-ndU-ndta}
    \\
    \cstf_{\Uan\ta}
    &\toprefb{eq:pU-pta-ndU-ndta}{eq:cstf-sstf}{= } 
    \tao^{-2}(\tao\pU+\moacs 1)\tao\pta\tao^{1-2\mv}\sstf
    =
    \tao^{-2\mv-1}\tao^{2\mv-1}(\tao\pU+\moacs 1)\tao^{1-2\mv}\tao^{2\mv-1}\tao\pta\tao^{1-2\mv}\sstf
    \\&\toprefb{eq:mdU}{eq:mdta}{=}
    \tao^{-2\mv-1}(\mdU+\moacs 1)\mdta\sstf
    \myeqlab{eq:cstf-ta-Uan-sstf}
    \\
    (\mdU+\moacs 1)\mdta\sstfz &\toprefb{eq:mdU}{eq:mdta-sstfz}{= } [(\tao\pU+\moacs2\moacs\mv-\moacs1)+\moacs 1](-1) 
    =  -2\mv \myeqlab{eq:vorfbase} 
\end{alignat*} 

\subsection{Nonlinear expressions}

\begin{alignat*}{5} \lrv_\Uan
    &\topref{eq:lrv-x}{\eq } 
    \frac{\cstf_{\Uan\ta}}{2\cstf_\ta}
    \toprefb{eq:cstf-ta-Uan-sstf}{eq:cstf-ta-sstf}{\eq }
    \frac{\tao^{-2\mv-1}(\mdU+\moacs 1)\mdta\sstf}{2\tao^{-2\mv}\mdta\sstf}
    \eq 
    \tao^{-1} \frac{(\mdU+\moacs 1)\mdta\sstf}{2\mdta\sstf}
    \myeqlab{eq:pU-lrv-scaled}
    \\
    \lrv_\uan
    &\topref{eq:lrv-x}{\eq } 
    \frac{\cstf_{\uan\ta}}{2\cstf_\ta}
    \toprefb{eq:cstf-ta-uan-sstf}{eq:cstf-ta-sstf}{\eq }
    \frac{\tao^{-2\mv}\mdu\mdta\sstf}{2\tao^{-2\mv}\mdta\sstf}
    \eq 
    \frac{\mdu\mdta\sstf}{2\mdta\sstf}
    \myeqlab{eq:pu-lrv-scaled}
\end{alignat*} 

\begin{alignat*}{5} \pT \cstf
    &\topref{eq:pTpa-ptapu}{\eq } 
    \lrv_\Uan^{-1}\pU \cstf
    \toprefb{eq:pU-lrv-scaled}{eq:cstf-U-sstf}{\eq } 
    \tao\moac{\frac{2\mdta\sstf}{(\mdU+\moacs 1)\mdta\sstf}} \tao^{-2\mv}\mdU\sstf
    \eq
    \tao^{1-2\mv}~\subeq{\moac{\frac{2\mdta\sstf}{(\mdU+\moacs 1)\mdta\sstf}} \mdU}{\eqdef\mdT}\ \sstf
    \myeqlab{eq:mdT}
    \\
    \zeroaccl{\mdT\sstf}
    &\eq  
    \frac{2\mdta\sstfz}{(\mdU+\moacs 1)\mdta\sstfz}\mdU\sstfz
    \toprefc{eq:mdta-sstfz}{eq:mdU-sstfz}{eq:vorfbase}{\eq}
    \frac{2(-1)}{-2\mv}\cdot 1
    \eq 
    \frac1\mv
    \myeqlab{eq:mdT-sstfz} 
    \\
    \liac{\mdT\sstf}
    &\topref{eq:mdT}{\eq}
    \liac{\frac{2\mdta\sstf~\mdU\sstf}{(\mdU+\moacs 1)\mdta\sstf}}
    \\&\eq
    \frac{2[\liac{\mdta\sstf}\zeroaccl{\mdU\sstf}+\zeroaccl{\mdta\sstf}\liac{\mdU\sstf}]}{\zeroaccl{(\mdU+\moacs 1)\mdta\sstf}}
    -
    \frac{2\zeroaccl{\mdta\sstf}\zeroaccl{\mdU\sstf}\liac{(\mdU+\moacs 1)\mdta\sstf}}{[\zeroaccl{(\mdU+\moacs 1)\mdta\sstf}]^2}
    \\&\toprefc{eq:mdta-sstfz}{eq:mdU-sstfz}{eq:vorfbase}{\eq}
    \frac{2[\mdta\lsstf\cdot1+(-1)\cdot\mdU\lsstf]}{-2\mv}
    -
    \frac{2\cdot(-1)\cdot1\cdot(\mdU+\moacs 1)\mdta\lsstf}{[-2\mv]^2}
    \\&\eq
    [ \frac{\mdU-\mdta}{\mv}
    +
    \frac{(\mdU+\moacs 1)\mdta}{2\mv^2} ] \lsstf
    \\&\eq
    \frac1{2\mv^2} [ (\mdU+\moacs 1)\mdta + 2\mv (\mdU-\mdta) ] \lsstf
    \eq
    \frac1{2\mv^2} [ 
    \mdU(\mdta+\moacs2\moacs\mv) - (2\mv-1) \mdta
    ] \lsstf 
    \myeqlab{eq:lin-stfTT}
\end{alignat*}

\begin{alignat*}{5} 
    \pa \cstf
    &\topref{eq:pTpa-ptapu}{\eq} 
    (\pu-\frac{\lrv_\uan}{\lrv_\Uan} \pU) \cstf
    \\&\topref{eq:lrv-x}{\eq }
    (\pu-\frac{\cstf_{\uan\ta}}{\cstf_{\Uan\ta}}\pU)\cstf
    \toprefd{eq:cstf-ta-uan-sstf}{eq:mdU}{eq:cstf-ta-Uan-sstf}{eq:cstf-sstf}{\eq} 
    (\tao^{1-2\mv}\mdu\tao^{2\mv-1} - \frac{\tao^{-2\mv}\mdu\mdta\sstf}{\tao^{-2\mv-1}(\mdU+\moacs 1)\mdta\sstf}\tao^{-2\mv}\mdU\tao^{2\mv-1})
    \tao^{1-2\mv}\sstf
    \\&\eq 
    \tao^{1-2\mv} ( \subeq{\mdu-\frac{\mdu\mdta\sstf}{(\mdU+\moacs 1)\mdta\sstf}\mdU}{\eqdef\mda} ) \sstf 
    \myeqlab{eq:mda}
\end{alignat*} 

\begin{alignat*}{5} \zeroaccl{\mda\sstf}
    &\topref{eq:mda}{\eq } 
    \mdu\sstfz - \frac{\mdu\mdta\sstfz}{(\mdU+\moacs 1)\mdta\sstfz}\mdU\sstfz
    \toprefd{eq:mdu-sstfz}{eq:mdu-mdta-sstfz}{eq:vorfbase}{eq:mdU-sstfz}{\eq} 
    0 - \frac{0}{-2\mv}\cdot1  \eq  0
    \myeqlab{eq:mda-sstfz} 
\end{alignat*} 

\begin{alignat*}{5}  
    \liac{\mda\sstf}
    &\topref{eq:mda}{\eq}
    \mdu\lsstf 
    - \frac{\liac{\mdu\mdta\sstf}\zeroaccl{\mdU\sstf}}{\zeroaccl{(\mdU+\moacs 1)\mdta\sstf}}
    - \frac{\zeroaccl{\mdu\mdta\sstf}\liac{\mdU\sstf}}{\zeroaccl{(\mdU+\moacs 1)\mdta\sstf}}
    + \frac{\zeroaccl{\mdu\mdta\sstf}\zeroaccl{\mdU\sstf}\liac{(\mdU+\moacs 1)\mdta\sstf}}{[\zeroaccl{(\mdU+\moacs 1)\mdta\sstf}]^2} 
    \\&\toprefd{eq:mdu-sstfz}{eq:mdu-mdta-sstfz}{eq:mdU-sstfz}{eq:vorfbase}{\eq}
    \mdu\lsstf 
    - \frac{\mdu\mdta\lsstf\cdot 1}{-2\mv}
    - \frac{0\cdot\mdU\lsstf}{-2\mv}
    + \frac{0\cdot 1\cdot(\mdU+\moacs 1)\mdta\lsstf}{(-2\mv)^2} 
    \\&\eq 
    \frac1{2\mv} \mdu(\mdta+\moacs2\moacs\mv)\lsstf
    \myeqlab{eq:lin-stf-a}
\end{alignat*}

\begin{alignat*}{5} 
    \gU
    \topref{eq:Delta-pav-new}{\eq} 
    (\pT-\lrv_\uan\pa)\cstf
    \toprefc{eq:mdT}{eq:pu-lrv-scaled}{eq:mda}{\eq}
    (\tao^{1-2\mv}\mdT-\frac{\mdu\mdta\sstf}{2\mdta\sstf}\tao^{1-2\mv}\mda)\sstf
    \eq 
    \tao^{1-2\mv}\subeq{[\mdT-\frac{\mdu\mdta\sstf}{2\mdta\sstf}\mda]\sstf}{\eqdef\sgU}
    \myeqlab{eq:sgU}
\end{alignat*} 

\begin{alignat*}{5}  
    \zeroaccl{\sgU} 
    &\topref{eq:sgU}{\eq } 
    \zeroaccl{\mdT\sstf}
    -
    \frac{\zeroaccl{\mdu\mdta\sstf}}{2\zeroaccl{\mdta\sstf}}
    \zeroaccl{\mda\sstf}
    \toprefb{eq:mdT-sstfz}{eq:mda-sstfz}{\eq}
    \frac1\mv 
    - 
    \frac{\zeroaccl{\mdu\mdta\sstf}}{2\zeroaccl{\mdta\sstf}}
    0
    \eq  
    \frac1\mv 
    \myeqlab{eq:gU-zero} 
\end{alignat*} 

\begin{alignat*}{5} 
    \liac{\sgU}
    &\topref{eq:sgU}{\eq } 
    \liac{ \mdT\sstf
        - \frac{\mdu\mdta\sstf}{2\mdta\sstf} 
        \mda\sstf }
    \topref{eq:lin-stfTT}{\eq }
    \liac{\mdT\sstf}
    -
    \liac{\frac{\mdu\mdta\sstf}{2\mdta\sstf}}
    \zeroaccl{\mda\sstf}
    -
    \frac{\zeroaccl{\mdu\mdta\sstf}}{2\zeroaccl{\mdta\sstf}}
    \liac{\mda\sstf}
    \\&\toprefc{eq:mdU-sstfz}{eq:vorfbase}{eq:mdu-mdta-sstfz}{\eq} 
    \liac{\mdT\sstf}
    -
    \liac{\frac{\mdu\mdta\sstf}{2\mdta\sstf}}
    \cdot 0
    -
    \frac{0}{2\zeroaccl{\mdta\sstf}}
    \liac{\mda\sstf}
    \\&\topref{eq:lin-stfTT}{\eq }
    \frac1{2\mv^2} [ 
    \mdU(\mdta+\moacs2\moacs\mv) - (2\mv-1) \mdta
    ] \lsstf 
    \myeqlab{eq:lin-gU}
\end{alignat*} 

\begin{alignat*}{5}
    \zeroaccl{\mdU \sgU} 
    \eq
    \mdU\zeroaccl{\sgU}
    \toprefb{eq:mdU}{eq:gU-zero}{\eq} 
    (\tao\pU+\moacs2\moacs\mv-\moacs1) \frac1\mv 
    \eq  
    \frac{2\mv-1}{\mv}
    \myeqlab{eq:mdU-gU-zero}
\end{alignat*} 

\begin{alignat*}{5}  
    \liac{\mdU\sgU}
    &\eq 
    \mdU\liac{\sgU}
    \topref{eq:lin-gU}{\eq}
    \frac1{2\mv^2} \mdU [ \mdU(\mdta+\moacs2\moacs\mv) - (2\mv-1) \mdta] \lsstf 
    \myeqlab{eq:lin-gUterm}
\end{alignat*} 

\begin{alignat*}{5} 
    \gu
    \topref{eq:Delta-pav-new}{\eq } 
    \lrv_\Uan\pa\cstf
    \toprefb{eq:pU-lrv-scaled}{eq:mda}{\eq }
    \tao^{-1}\frac{(\mdU+\moacs 1)\mdta\sstf}{2\mdta\sstf}\tao^{1-2\mv}\mda\sstf
    \eq 
    \tao^{-2\mv}\subeq{\frac{(\mdU+\moacs 1)\mdta\sstf}{2\mdta\sstf}\mda\sstf}{\eqdef\sgu}
    \myeqlab{eq:sgu}
\end{alignat*} 

\begin{alignat*}{5}
    \zeroaccl{\sgu}
    &\topref{eq:sgu}{\eq }
    \frac{(\mdU+\moacs 1)\mdta\sstfz}{2\mdta\sstfz} \mda\sstfz
    \topref{eq:mda-sstfz}{\eq}
    \frac{(\mdU+\moacs 1)\mdta\sstfz}{2\mdta\sstfz} 0
    \eq 
    0
    \myeqlab{eq:gu-zero} 
\end{alignat*} 

\begin{alignat*}{5} 
    \liac{\sgu}
    &\topref{eq:sgu}{\eq } 
    \liac{
        \frac{(\mdU+\moacs 1)\mdta\sstf}{2\mdta\sstf}
        \mda\sstf
    }
    \eq 
    \liac{ \frac{(\mdU+\moacs 1)\mdta\sstf}{2\mdta\sstf} }
    \zeroaccl{\mda\sstf}
    +
    \zeroaccl{ \frac{(\mdU+\moacs 1)\mdta\sstf}{2\mdta\sstf} }
    \liac{\mda\sstf}
    \\&\toprefd{eq:mda-sstfz}{eq:vorfbase}{eq:mdta-sstfz}{eq:lin-stf-a}{\eq}
    \liac{ \frac{(\mdU+\moacs 1)\mdta\sstf}{2\mdta\sstf} }
    0
    +
    \frac{-2\mv}{2(-1)}
    [\frac1{2\mv} \mdu(\mdta+\moacs2\moacs\mv)\lsstf]
    \eq  
    \frac12 \mdu(\mdta+\moacs2\moacs\mv)\lsstf
    \myeqlab{eq:lin-gu}
\end{alignat*} 

\begin{alignat*}{5}
    \zeroaccl{\mdu \sgu} 
    \eq
    \mdu\zeroaccl{\sgu}
    \topref{eq:gu-zero}{\eq } 
    \mdu 0 
    \eq  
    0
    \myeqlab{eq:mdu-gu-zero}
\end{alignat*} 

\begin{alignat*}{5} 
    \liac{\mdu\sgu}
    &\eq
    \mdu\liac{\sgu}
    \topref{eq:lin-gu}{\eq}
    \frac12 \mdu^2(\mdta+\moacs2\moacs\mv)\lsstf
    \myeqlab{eq:lin-guterm}
\end{alignat*}

\begin{alignat*}{5} 
    \cvorf
    \topref{eq:vort-stf-ta}{\eq}
    \cstf_\Uan^{-\frac1{2\mv}}
    \topref{eq:cstf-U-sstf}{\eq}
    (\tao^{-2\mv}\mdU\sstf)^{-\frac1{2\mv}} 
    \eq
    \tao \subeq{(\mdU\sstf)^{-\frac1{2\mv}} }{\eqdef\svorf}
    \myeqlab{eq:cvorf-svorf}
\end{alignat*} 

\begin{alignat*}{5} 
    \zeroacc\svorf
    &\topref{eq:cvorf-svorf}{\eq}
    (\mdU\sstfz)^{-\frac1{2\mv}}
    \topref{eq:mdU-sstfz}{\eq}
    (1)^{-\frac1{2\mv}}
    \eq
    1
    \myeqlab{eq:svorfz}
\end{alignat*} 

\begin{alignat*}{5} 
    \liac{\svorf}
    &\topref{eq:cvorf-svorf}{\eq}
    \liac{(\mdU\sstf)^{-\frac1{2\mv}}}
    \eq
    -\frac1{2\mv} \zeroaccl{(\mdU\sstf)^{-\frac1{2\mv}-1}} \mdU\lsstf
    \topref{eq:mdU-sstfz}{\eq}
    -\frac1{2\mv} (1)^{-\frac1{2\mv}-1} \mdU\lsstf
    \\&\eq
    -\frac1{2\mv} \mdU\lsstf 
    \myeqlab{eq:svorf-lin}
\end{alignat*} 

\begin{alignat*}{5} 
    \cvort
    \topref{eq:vort-stf-ta}{\eq}
    \cvorf \Vort
    \topref{eq:cvorf-svorf}{\eq}
    \tao \subeq{\svorf \Vort}{\eqdef\svort}
    \myeqlab{eq:cvort-svort}
\end{alignat*} 

\begin{alignat*}{5} \svortz
    &\topref{eq:cvort-svort}{\eq}
    \zeroaccl{\svorf\sVort}
    \toprefb{eq:svorfz}{eq:Vortz}{\eq}
    1\cdot\frac{2\mv-1}{\mv}
    =
    \frac{2\mv-1}{\mv}
    \myeqlab{eq:svort-zero}
\end{alignat*} 

\begin{alignat*}{5} \liac{\svort}
    &\topref{eq:cvort-svort}{\eq}
    \liac{\svorf\sVort}
    \eq
    \liac{\svorf}\sVortz+\zeroacc\svorf\Vortl
    \toprefc{eq:svorf-lin}{eq:Vortz}{eq:svorfz}{\eq}
    -\frac1{2\mv}\mdU\lsstf\frac{2\mv-1}{\mv}+1\cdot\sVortl
    \\&=
    \frac{1-2\mv}{2\mv^2}\mdU\lsstf+\sVortl
    \myeqlab{eq:svort-lin}
\end{alignat*} 

\begin{alignat*}{5}
    \vortterm
    \topref{eq:Delta-pav-new}{\eq}
    - \frac1{2\mv} \cstf_{\Uan\ta} \cdot \cvort
    \toprefb{eq:cstf-ta-Uan-sstf}{eq:cvort-svort}{\eq }
    - \frac1{2\mv} \tao^{-2\mv-1} (\mdU+\moacs 1)\mdta\sstf \cdot \tao \svort
    \eq 
    \tao^{-2\mv} \subeq{[ - \frac1{2\mv} (\mdU+\moacs 1)\mdta\sstf \cdot \svort]}{\eqdef\svortterm}
    \myeqlab{eq:vortterm-new}
\end{alignat*}

\begin{alignat*}{5}
    \zeroacc\svortterm
    &\topref{eq:vortterm-new}{\eq }
    -\frac1{2\mv} (\mdU+\moacs 1)\mdta\sstfz \cdot \svortz
    \toprefb{eq:vorfbase}{eq:svort-zero}{\eq}
    -\frac1{2\mv} (-2\mv) \frac{2\mv-1}{\mv} 
    \eq 
    \frac{2\mv-1}{\mv} 
    \myeqlab{eq:vortterm-zero} 
\end{alignat*}

\begin{alignat*}{5}  
    \liac{\svortterm}
    &\topref{eq:vortterm-new}{\eq}
    \liac{-\frac1{2\mv}(\mdU+\moacs 1)\mdta\sstf\cdot\svort}
    \eq
    -\frac1{2\mv}[(\mdU+\moacs 1)\mdta\lsstf\cdot\svortz+(\mdU+\moacs 1)\mdta\sstfz\cdot\liac{\svort}]
    \\&\toprefc{eq:vorfbase}{eq:svort-zero}{eq:svort-lin}{\eq}
    -\frac1{2\mv}[(\mdU+\moacs 1)\mdta\lsstf\frac{2\mv-1}{\mv}+(-2\mv)(\frac{1-2\mv}{2\mv^2}\pU\lsstf+\sVortl)]
    \\&\eq
    \frac{1-2\mv}{2\mv^2}
    [ \mdU
    +
    \mdU\mdta
    +
    \mdta ] \lsstf
    +
    \Vortl 
    \myeqlab{eq:lin-vortterm}
\end{alignat*} 

\subsection{Other quantities}

\begin{alignat*}{5} 
    \crad
    \topref{eq:rad-lrv}{=}
    e^{\lrv}
    \topref{eq:lrv-cstf}{=}
    (\frac{\pta\cstf}{-\mv})^{\frac12}
    \topref{eq:cstf-ta-sstf}{=}
    \tao^{-\mv}~\subeq{(\frac{\mdta\sstf}{-\mv})^{\frac12}}{\eqdef\srad}
    \myeqlab{eq:crad-srad}
\end{alignat*} 

\begin{alignat*}{5} \sradz 
    &\topref{eq:crad-srad}{\eq} 
    (\frac{\mdta\sstfz}{-\mv})^{\frac12}
    \topref{eq:mdta-sstfz}{=} 
    (\frac{-1}{-\mv})^{\frac12} 
    = 
    \mv^{-\frac12} \myeqlab{eq:sradz} \end{alignat*} 

\begin{alignat*}{5} \liac{\srad} 
    &\topref{eq:crad-srad}{\eq} 
    \liac{ (\frac{\mdta\sstf}{-\mv})^{\frac12} }
    \eq
    \frac12 \mv^{-\frac12} \zeroaccl{(-\mdta\sstf)^{-\frac12}} \mdta\lsstf
    \topref{eq:mdta-sstfz}{\eq}
    \frac12 \mv^{-\frac12} \mdta\lsstf
    \myeqlab{eq:srad-lin}
\end{alignat*} 

\subsection{Overall map}

\begin{alignat*}{5} 0 
    &\topref{eq:Delta-pav-new}{\eq}
    2\mv^2 \tao^{2\mv} [ -\pU\gU - \pu\gu + \vortterm ] 
    \toprefc{eq:sgU}{eq:sgu}{eq:vortterm-new}{\eq } 
    2\mv^2 \tao^{2\mv} [ -\pU\tao^{1-2\mv}\sgU - \pu\tao^{-2\mv}\sgu + \tao^{-2\mv}\svortterm]
    \\&\topref{eq:mdU}{\eq } 
    2\mv^2 (-\mdU\sgU-\pu\sgu+\svortterm)
    \eqdef\Fmap
    \myeqlab{eq:curl-scaled-new}
\end{alignat*} 

The following step proves that $\stfz$, the background solution, is actually a solution of the nonlinear PDE
\myeqref{eq:Delta-pav-new}.
\begin{alignat*}{5}
    \zeroacc{\Fmap}
    \eq
    2\mv^2 [- \zeroaccl{\mdU \sgU}
    - \zeroaccl{\mdu \sgu}
    + \zeroacc{\svortterm}]
    \toprefc{eq:mdU-gU-zero}{eq:mdu-gu-zero}{eq:vortterm-zero}{\eq}
    2\mv^2 [ - \frac{2\mv-1}{\mv} - 0 + \frac{2\mv-1}{\mv} ]
    \eq  
    0
    \myeqlab{eq:background-solves}
\end{alignat*}
\begin{alignat*}{5} &
    \liac{\Fmap}
    \topref{eq:curl-scaled-new}{\eq}
    2\mv^2 [- \liac{ \mdU \sgU} - \liac{\mdu \sgu} + \liac{\svortterm } ]
    \\\toprefc{eq:lin-gUterm}{eq:lin-guterm}{eq:lin-vortterm}{\eq}&
    2\mv^2 \Big[ 
    -\Big( 
    \frac1{2\mv^2} \mdU [ \mdU(\mdta+\moacs2\moacs\mv) - (2\mv-1) \mdta ] \lsstf
    \Big)
    -
    \Big( \frac12 \mdu^2(\mdta+\moacs2\moacs\mv) ] \lsstf \Big)
    \\&\qquad+
    \Big( 
    \frac{1-2\mv}{2\mv^2} 
    [ \mdU
    +
    \mdU\mdta
    +
    \mdta] \lsstf 
    + \sVortl \Big)
    \Big]
    \\\eq&
    \Big(
    -
    \mdU^2(\mdta+\moacs2\moacs\mv) 
    +
    (2\mv-1) \mdU\mdta 
    -
    \mv^2\mdu^2(\mdta+\moacs2\moacs\mv) 
    -
    (2\mv-1) [
    \mdU
    +
    \mdU\mdta
    +
    \mdta
    ]
    \Big) 
    \lsstf 
    \\&\qquad+ 
    2\mv^2
    \sVortl 
    \\\eq&
    \subeq{
        \big(
        (-\mdU^2-\mv^2\mdu^2) (\mdta+\moacs2\moacs\mv) 
        -
        (2\mv-1) (\mdU + \mdta)
        \big) 
    }{=\frac{\partial F}{\partial\sstf}(\sstfz,\sVortz)\eqdef\Lop}
    \lsstf 
    + 
    \subeq{
        2\mv^2 
    }{=\frac{\partial F}{\partial\sVort}(\sstfz,\sVortz)}
    \sVortl 
    \myeqlab{eq:linearization}
\end{alignat*}

\subsection{Fourier transform}

\mylabel{section:foutr}%

While the background solution $\stfz$ has algebraic spiral pseudo-streamlines, any perturbed solution will
have ``wobbly'' spirals, with deformations such as ellipsoidal stretching.
Such perturbations tend to be pronounced at ``integer frequencies'', since non-integer frequency perturbations favor
self-intersection of the spiral. Compare, for example, 
\begin{alignat*}{5} z(\ta)=\ta^{-\mv}e^{i\ta}(1+\delta~\ta^{-\alpha} e^{i\ta\etav }) \myeqlab{eq:selfint} \end{alignat*} 
in Figure \myref{fig:self-intersection} for integer and non-integer $\etav$ ($\delta=0$ corresponds to the 
pseudo-streamlines of the background solution).
\if\dofigures%
\begin{figure}
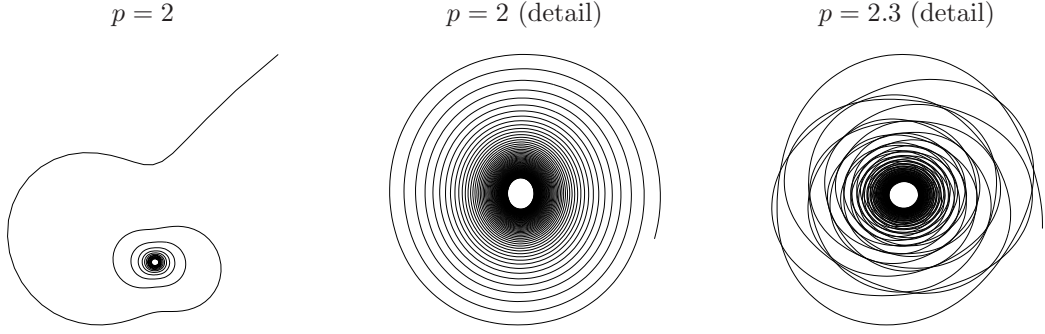

    \centerline{\parbox{2in}{\input{selfintersectionlarge.pstex}}\parbox{2in}{\input{selfintersection2.pstex}}\parbox{2in}{\input{selfintersection2half.pstex}}}
    \caption{Perturbed spirals; $\mv=1,\alpha=0.5,\delta=0.7$. (The spiral center is not drawn, leaving a white spot in the middle.) 
        Integer frequency $\etav$: the perturbation is barely noticeable near the center. Non-integer frequency $\etav$: 
        physically unreasonable self-intersection. These will occur for any $\alpha<1$ and $\delta>0$, if sufficiently large $\ta$ are considered.
}
    \mylabel{fig:self-intersection}
\end{figure}
\fi%

On the other hand the $\vv$ field induced by non-integer frequency perturbations will be weaker because the phase-shift between adjacent spiral turns causes their contributions to cancel each other. But for near-integer frequencies there is \defm{resonance}; the contributions
add up.

Hence we are motivated to take the Fourier transform. Technically this has the convenient
effect of reducing the equation to a family of \emph{ordinary} differential
equations, with the derivatives stemming from the (unavoidable) presence
of $\ta$ (hidden in $\mdta,\mdU$, see \myeqref{eq:mdU} and \myeqref{eq:mdta}) 
in the scaled $\ta,\uan$-space equation \myeqref{eq:linearization}).

We use $(\etav,\nv)\in\R\times\Z$ as the Fourier variables for $\ta,\uan\in\R$
(for technical convenience we consider functions defined for all $\ta\in\R$, although their physical relevance is only for $\ta>0$):
\begin{alignat*}{5} 
    (\mdu\lsstf)^\wedge 
    &\eq  
    i\nv\lsstff \myeqlab{eq:mdu-fou} 
\end{alignat*} 
\begin{alignat*}{5}
    (\mdta\lsstf)^\wedge 
    &\topref{eq:mdta}{\eq}
    [(\tao\pta+\moacs1-\moacs2\moacs\mv)\lsstf]^\wedge 
    \eq  
    (i\peta i\etao+\moacs1-\moacs2\moacs\mv)\lsstff 
    \eq  
    (-\peta\etao+\moacs1-\moacs2\moacs\mv)\lsstff 
    \\&\eq  
    -(\subeq{\etao\peta}{\eqdef \Bo}+\moacs2\moacs\mv)\lsstff 
    \myeqlab{eq:mdta-fou-a} 
\end{alignat*} 
\begin{alignat*}{5}
    \big((\mdta+\moacs2\moacs\mv)\lsstf\big)^\wedge 
    &\topref{eq:mdta-fou-a}{\eq} 
    -\etao\peta\lsstff  
    \topref{eq:mdta-fou-a}{\eq} 
    -\Bo\lsstff  
    \myeqlab{eq:mdta-fou}
\end{alignat*} 
\begin{alignat*}{5}
    (\mdU\lsstf)^\wedge 
    &\topref{eq:mdU}{\eq}
    [(\tao\pU+\moacs2\moacs\mv-\moacs1)\lsstf]^\wedge 
    \\&\topref{eq:pU-def}{\eq}
    \big([\tao(\pu-\pta)+\moacs2\moacs\mv-\moacs1]\lsstf\big)^\wedge 
    \\&\eq  
    [i\peta(i\nvo-i\etao)+\moacs2\moacs\mv-\moacs1]\lsstff 
    \\&\eq  
    [\peta(\etao-\nvo)+\moacs2\moacs\mv-\moacs1]\lsstff 
    \eq  
    \big(\subeq{(\etao-\nvo)\peta}{\eqdef \Ao}+\moacs2\moacs\mv\big)\lsstff  
    \myeqlab{eq:mdU-fou} 
\end{alignat*} 
\begin{alignat*}{5}
    \Big(\frac{\partial\Fmap}{\partial\sstf}[\sstfz,\Vortz]\lsstf\Big)^\wedge
    &\topref{eq:linearization}{\eq}
    \Big(\big(
    [-\mdU^2-\mv^2\mdu^2] (\mdta+\moacs2\moacs\mv) 
    -
    (2\mv-1) [\mdU + \mdta]
    \big) 
    \lsstf \Big)^\wedge
    \\&\toprefd{eq:mdU-fou}{eq:mdu-fou}{eq:mdta-fou}{eq:mdta-fou-a}{\eq}
    \big(
    [-(\Ao+\moacs2\moacs\mv)^2
    - \moacs\mv^2(i\moacs\nv)^2]
    (-\Bo) 
    - (2\mv-1) [(\Ao+\moacs2\moacs\mv)-(\Bo+\moacs2\moacs\mv)]
    \big)
    \lsstff 
    \\&=
    \big(
    [(\Ao+\moacs2\moacs\mv)^2
    - (\moacs\mv\nvo)^2]
    \Bo
    - (\moacs2\moacs\mv-\moacs1) [\Ao-\Bo]
    \big)
    \lsstff 
    \\&=
    \big(
    \subeq{
        \supeq{[\Ao+\supeq{(\moacs2-\moacs\nv)\moacs\mv}{\eqdef\mmen}]
            [\Ao+\supeq{(\moacs2+\moacs\nv)\moacs\mv}{\eqdef\mmep}]
            \Bo}{\eqdef\Rop}
        -
        \supeq{(\moacs2\moacs\mv-\moacs1) (\Ao-\Bo)}{\eqdef\Eop}
    }{\eqdef\Lopn}
    \big)
    \lsstff 
    \myeqlab{eq:lstff-new} 
\end{alignat*}  

\section{Integral operators}%
\mylabel{section:intops}%

$\Aopn$ and $\Bo$ in \myeqref{eq:mdU-fou} and \myeqref{eq:mdta-fou-a} are of the type $\lops$ for $\pole,\exs\in\R$, so we seek inverses of such operators now.
We discuss only the non-pathological case $\exs\neq-1$
and consider only $u\in\Mh$. (It is possible to discuss $\Ms$, but not needed for our purposes.)

\subsection{Uniqueness}

Consider $u\in\Ms(\R)$ so that
\begin{alignat*}{5} \lops u = 0\quad. \end{alignat*} 
Then
\begin{alignat*}{5} \peta u = (\etao-\polo)^{-1}\exso u\quad, \end{alignat*} 
so $u\in\Ms(\R)$ yields $u\in\BV_\loc(\R\backslash\{\pole\})$ which in turn yields higher regularity until we have
$u\in\Cinf(\R\backslash\{\pole\})$. Thus $u$ is a classical solution,  
hence a multiple of $|\etav-\pole|^\exs$, on each side of $\pole$. 
This is not integrable for any $\exs\in\R$ (for $\exs\geq-1$ at $\etav=\pm\infty$,
for $\exs\leq-1$ at $\etav=\pole\pm$). Hence $u=0$ on each side of $\pole$, so $u$ is a measure supported in $\pole$.
Thus $u=\alpha\delta(\cdot-\pole)$ for some $\alpha\in\CC$. But 
\begin{alignat*}{5} &0 = \lops u = \alpha \lops \delta(\etav-\pole) 
    = \alpha [ \peta\subeq{(\etao-\polo)\delta(\etav-\pole)}{=0}-\subeq{(\exso+1)}{\neq 0}\delta(\etav-\pole) ]
    \\&\quad\impl\quad \alpha = 0 
\end{alignat*} 
Hence $u=0$.

(Uniqueness need not hold in larger $u$ classes; e.g.\ $\etao\peta u=0$ is solved by $u=\cf_{\boi0\infty}$.)

\subsection{Existence}

For $\xiv\in\R$ consider 
\begin{alignat*}{5} \lops u=\delta(\cdot-\xiv) \end{alignat*}
For $\xiv=\pole$ the solution is simply
\begin{alignat*}{5} u = -(\exs+1)^{-1}\delta(\cdot-\xiv) \end{alignat*} 
since
\begin{alignat*}{5} \lops \delta(\cdot-\pole)
    &= \lopsd \delta(\cdot-\pole)
    \\&= \peta\subeq{(\etao-\polo)\delta(\cdot-\pole)}{=0}-(\exs+1)\delta(\cdot-\pole)
    = -(\exs+1)\delta(\cdot-\pole) \end{alignat*} 
Consider $\xiv\neq\pole$. We begin with $\fff\in\LoneR$.

For $\etav$ on one side of $\pole$,
\begin{alignat*}{5} \fff
    &= 
    \lops u
    =
    \moac{|\etav-\pole|^{\exs+1}} \peta \moac{\sign(\etav-\pole)|\etav-\pole|^{-\exs}} u
    \\\quad\impl\quad
    u(\etav) &= \sign(\etav-\pole) |\etav-\pole|^{\exs} \Big( C + \int_{c}^{\etav-\pole} |\xiv-\pole|^{-\exs-1} \fff(\xiv)d\xiv \Big) \end{alignat*} 
We combine the two sides and choose $C,c$ appropriately:
\begin{alignat*}{5} u(\etav) 
    &\defeq
    \sign(\etav-\pole)|\etav-\pole|^{\exs}\int_{\sign(\etav-\pole)\infty^{\exs+1}}^{\etav-\pole} |\xiv-\pole|^{-\exs-1}\fff(\xiv)d(\xiv-\pole) 
    \myeqlab{eq:uint}
    \\&= \int_{\infty^{\exs+1}}^1 \xv^{-\exs-1} \fff(\xv\etav) d\xv \end{alignat*} 
By choice of lower integral boundary, $|\xiv-\pole|^{-\exs-1}$ in \myeqref{eq:uint} is bounded,
while $\fff$ is integrable, so $u$ is well-defined for $\etav\neq\pole$.
\begin{alignat*}{5} \norma{u}{\LoneR}
    &= \int_\R|u(\etav)|d\etav 
    = \int_\R \Big| \int_{\infty^{\exs+1}}^1 \xv^{-\exs-1} \fff(\xv\etav) d\xv \Big| d\etav
    \\&\leq \int_\R \int_{\infty^{\exs+1}}^1 \xv^{-\exs-1} |\fff(\xv\etav)| d\xv~d\etav
    = \int_{\infty^{\exs+1}}^1 \xv^{-\exs-1} \int_\R |\fff(\xv\etav)| d\etav~d\xv
    \\&= \int_{\infty^{\exs+1}}^1 \xv^{-\exs-1} \xv^{-1} \int_\R |\fff(\xiv)| d\xiv~d\xv
    \\&= |\exs+1|^{-1}\norma{\fff}{\LoneR}\quad.
\end{alignat*} 
Fundamental theorem of calculus:
\begin{alignat*}{5} 
    \lops u(\etav) 
    \topref{eq:uint}{=}&
    (\etav-\pole) \Big( \exs|\etav-\pole|^{\exs-1}
    \int_{\sign(\etav-\pole)\infty^{\exs+1}}^{\etav-\pole} |\xiv-\pole|^{-\exs-1} \fff(\xiv) d(\xiv-\pole)
    \\&\qquad+ \sign(\etav-\pole)|\etav-\pole|^\exs \big[ |\xiv-\pole|^{-\exs-1} \fff(\xiv) \big]_{|\xiv-\pole=\etav-\pole} \Big)
    \\&
    -\exs 
    \sign(\etav-\pole) |\etav-\pole|^{\exs}\int_{\sign(\etav-\pole)\infty^{\exs+1}}^{\etav-\pole} |\xiv-\pole|^{-\exs-1}\fff(\xiv)d(\xiv-\pole) 
    =
    \fff(\etav)
\end{alignat*} 
Hence 
\begin{alignat*}{5} \lops u=\fff \quad\text{on $\R\backslash\{\pole\}$} \myeqlab{eq:lopsu-nopole} \end{alignat*} 
in classical and in distributional sense; it remains to discuss $\etav=\pole$.

Let $\phi\in\Cinf_c(\R)$ be some test function. We prepare to cut out a neighbourhood of $\pole$ from its support
by choosing a $\theta_1\in\Cinf(\R)$ with value $0$ on $\cli{-1}{1}$ and value $1$ outside $\cli{-2}{2}$.
Define $\theta_\eps(\cdot):=\theta_1(\frac{\cdot-\pole}{\eps})$. Then $\theta_\eps(\etav)=1$ for $\etav\not\in\pole+\cli{-2\eps}{2\eps}$, 
$\peta\theta_\eps(\etav)=0$ for $\etav\not\in\pole+\cli{-2\eps}{2\eps}\setdiff\cli{-\eps}{\eps}$, 
and $\norma{\peta\theta_\eps}{\infty}\lesssim\eps^{-1}$ as $\eps\searrow 0$, so $\norma{(\etao-\polo)\peta\theta_\eps}{\infty}\lesssim 1$. 
\begin{alignat*}{5} 
    \langle\phi,\fff\rangle
    &=
    \int_\R\phi\fff
    \overset{\eps\searrow 0}{\leftarrow}
    \int_\R \theta_\eps \phi \fff 
    \topref{eq:lopsu-nopole}{=} 
    \int_\R \theta_\eps \phi \lops u 
    \\&=
    \int_\R u \big(-\peta[\theta_\eps(\etao-\polo)\phi]\big) - \exs u\theta_\eps\phi
    \\&=
    \int_\R u \theta_\eps (-\peta(\etao-\polo)-\exso) \phi
    +
    \int_\R \subeq{-u \phi}{\in\LoneR}\,\subeq{(\etao-\polo)\peta\theta_\eps}{\anorma{\infty}\lesssim 1,\supp\subset\pole+\cli{-2\eps}{2\eps}}
    \\&\overset{\eps\searrow0}{\rightarrow}
    \int_\R u (-\peta(\etao-\polo)-\exso) \phi
    =
    \big\langle \lops u , \phi \big\rangle
\end{alignat*} 
Hence $u$ solves $\lops u=\fff$ on $\R$ in the distributional sense.

For $\fff=\delta(\etav-\xiv)$ for $\xiv\neq\pole$ we repeat the analysis with obvious modifications.

Altogether we have a unique solution $u=\Keps\fff\in\Mh$ for any $\fff\in\Mh$, with operator norm estimate
\begin{alignat*}{5} \normM{\Keps} \leq |\exs+1|^{-1} \myeqlab{eq:Keps-normM} \end{alignat*} 

\section{Inversion for each $\nv$ separately} 

\mylabel{section:inversion}

\subsection{$\Ao,\Bo$ commutators}

We can commute $\Ao+\mep$, $\Ao+\men$ freely, but not necessarily $\Ao,\Bo$. 
\begin{alignat*}{5}
    \Ao\Bo
    &\toprefb{eq:mdU-fou}{eq:mdta-fou-a}{\eq}
    (\etao-\nvo)\peta \etao\peta 
    \\&\eq  (\etao-\nvo)(\etao\peta+\moacs1)\peta 
    \\&\eq \etao(\etao-\nvo)\peta\peta+(\etao-\nvo)\peta
    \\&\eq \etao\big(\peta(\etao-\nvo)-\moacs1\big)\peta+(\etao-\nvo)\peta
    \\&\eq  \etao\peta(\etao-\nvo)\peta - \etao\peta + (\etao-\nvo)\peta 
    \\&\eq  \Bo\Ao + \Ao - \Bo  
    \myeqlab{eq:ABBA}
    \\\eqv\quad 
    \Bo \Ao  &\eq  \Ao \Bo  + \Bo  - \Ao  
    \myeqlab{eq:BAAB}
\end{alignat*}
so
\begin{alignat*}{5}
    (\Ao +\mgr)(\Bo +\mgt) 
    &\eq  \Ao \Bo +\mgt\Ao +\mgr\Bo +\mgr\mgt  
    \\&\topref{eq:ABBA}{\eq } \big( \Bo \Ao  + \Ao  - \Bo\big)  + \mgt\Ao  + \mgr\Bo  + \mgr\mgt  
    \\&\eq  \Bo \Ao  + (\mgt+\moacs1)\Ao  + (\mgr-\moacs1)\Bo  + \mgr\mgt 
    \\&\eq  [\Bo +(\mgt+\moacs1)][\Ao +(\mgr-\moacs1)] - (\mgt+\moacs1)(\mgr-\moacs1) + \mgr\mgt 
    \\&\eq  [\Bo +\mgt+\moacs1][\Ao +\mgr-\moacs1] - \mgr + \mgt + \moacs1
    \myeqlab{eq:arbs}
\end{alignat*}
and, since $\Bo,\Ao$ can be interchanged in \myeqref{eq:ABBA}, 
\begin{alignat*}{5} 
    (\Bo +\mgt)(\Ao +\mgr)  
    \eq  [\Ao +\mgr+\moacs1][\Bo +\mgt-\moacs1] - \mgt + \mgr + \moacs1
    \myeqlab{eq:bsar}
\end{alignat*} 
\begin{alignat*}{5} &(\Ao +\mgr)(\Ao +\mgs)(\Bo +\mgt)
    \\&\topref{eq:arbs}\eq  (\Ao +\mgr)[(\Bo +\mgt+\moacs1)(\Ao +\mgs-\moacs1)-\mgs+\mgt+\moacs1]
    \\&\eq  (\Ao +\mgr)(\Bo +\mgt+\moacs1)(\Ao +\mgs-\moacs1)+(\mgt-\mgs+\moacs1)(\Ao +\mgr)
    \\&\topref{eq:arbs}{\eq}  [(\Bo +\mgt+\moacs2)(\Ao +\mgr-\moacs1)-\mgr+(\mgt+\moacs1)+\moacs1](\Ao +\mgs-\moacs1)+(\mgt-\mgs+\moacs1)(\Ao +\mgr)
    \\&\eq  (\Bo +\mgt+\moacs2)(\Ao +\mgr-\moacs1)(\Ao +\mgs-\moacs1)+(\mgt-\mgr+\moacs2)(\Ao +\mgs-\moacs1)+(\mgt-\mgs+\moacs1)(\Ao +\mgr)
    \\&\eq  (\Bo +\mgt+\moacs2)(\Ao +\mgr-\moacs1)(\Ao +\mgs-\moacs1)+(\moacs2\mgt-\mgr-\mgs+\moacs3)\Ao +[(\mgt-\mgr+\moacs2)(\mgs-\moacs1)+(\mgt-\mgs+\moacs1)\mgr]
    \\&\eq  (\Bo +\mgt+\moacs2)(\Ao +\mgr-\moacs1)(\Ao +\mgs-\moacs1)+(\moacs2\mgt-\mgr-\mgs+\moacs3)\Ao +[\mgt(\mgr+\mgs-\moacs1)-\moacs2(\mgr-\moacs1)(\mgs-\moacs1)]
    \myeqlab{eq:arasbt}
\end{alignat*} 

\subsection{$\Rop$ inversion and regularity estimates (Approach from non-commuted form)}

\subsubsection{Invertibility}

Consider the operator $\Rop=(\Aop)(\Aon)\Bo$ from \myeqref{eq:lstff-new}. 
As shown in Section \myref{section:intops}, 
\begin{alignat*}{5} \Aopn:(\Aopn)^{-1}[\Mh]\rightarrow\Mh \end{alignat*} 
(where the domain carries the induced norm) is an isomorphism
if $\exs=-\mepn\neq-1$. 
$\mep\topref{eq:lstff-new}{=}(2+\nv)\mv\eq1$ holds if and only if $\nv\neq -2$ and $\mv\eq\frac1{2+\nv}$. 
For $\nv<-2$ this yields $\mv<0$ which is outside the range $\boi{\frac12}\infty$ we are considering (see Theorem \myref{th:mainresult}).
For $\nv\geq 0$ we have $\mv\leq\frac12$ as well. 
For $\nv\eq -1$ there is a single solution $\mv\eq 1$. Hence $\Aop$ is invertible if $\nv\neq-1$ or $\mv\neq 1$. 

Analogously, $\Aon$ is invertible if $\nv\neq+1$ or $\mv\neq 1$. $\Bo$ is always invertible.

We will only consider functions whose Fourier transforms are zero for $\nv\not\in\Nv\Z$, corresponding to 
inverse transforms that are $2\pi/\Nv$-periodic in $\uan$.
We restrict $\Nv\geq 2$ 
so that $\Aopn$ is always invertible 
($\nv\in\Nv\Z$ means $|\nv|\neq 1$). 
$\Rop$ (defined on $\Rop^{-1}[\Mh]$ into $\Mh$) has 
inverse 
\begin{alignat*}{5} \Rop^{-1} \eq \Bo^{-1}(\Aon)^{-1}(\Aop)^{-1}  \quad. \myeqlab{eq:Rop-inverse} \end{alignat*} 

\subsubsection{Inverse norms}

For now let $\lesssim$ and $\sim$ be for the limit $\nv\rightarrow\infty$. 
$\gr,\gs,\gt$ will represent functions of $\nv$ so that 
\begin{alignat*}{5} \gr,\gs\lesssim\nv \quad,\quad \gt\lesssim 1 \myeqlab{eq:grgsgt} \end{alignat*}
All function norms are for $\Mh$ and all operator norms are for $[\Mh]$.
When using $\fff=\Bo^{-1}\Bo\fff$ etc., it is important to verify that $\Bo\fff$ is still in $\Mh$
so that $\Bo^{-1}\Bo\fff$ is well-defined; to avoid clutter we do not point this out explicitly every time. ($\fff=\Bo\Bo^{-1}\fff$
only requires $\fff\in\Mh$, on the other hand.)
We use these notations for the rest of Section \myref{section:inversion}. 
\renewcommand{\normm}[1]{\norm{#1}}
\renewcommand{\normM}[1]{\norm{#1}}

\begin{alignat*}{5} |\mepn| \topref{eq:lstff-new}{\sim} \asy\nv \myeqlab{eq:mepn-bound} \end{alignat*}
\begin{alignat*}{5} \normM{(\Aopn)^{-1}} \topref{eq:mdU-fou}{=} \normM{((\etao-\nvo)\peta-(-\mepn))^{-1}} 
    \topref{eq:Keps-normM}{\leq} |1+(-\mepn)|^{-1} \topref{eq:mepn-bound}{\lesssim} \asy\nv^{-1} \myeqlab{eq:Aopn-norm} \end{alignat*} 
\begin{alignat*}{5} \normM{(\Ao+\nv)^{-1}}
    \topref{eq:mdU-fou}{=} \normM{((\etao-\nvo)\peta-(-\nv))^{-1}}
    \topref{eq:Keps-normM}{\leq} |1+(-\nv)|^{-1} 
    \topref{eq:mepn-bound}{\lesssim} 
    \asy\nv^{-1} 
    \myeqlab{eq:Aoplusn-norm} 
\end{alignat*} 
and
\begin{alignat*}{5} \normM{\Bo^{-1}} &\topref{eq:mdta-fou}{=} \normM{[\etao\peta-0]^{-1}} \topref{eq:Keps-normM}{\leq} |1+0|^{-1}\eq 1 \myeqlab{eq:B-norm} \\
    \normM{(\Bo+\moacs2)^{-1}} &\topref{eq:mdta-fou}{=} \normM{[\etao\peta-(-\moacs2)]^{-1}} \topref{eq:Keps-normM}{\leq} |1+(-2)|^{-1}\eq 1 \myeqlab{eq:Btwo-norm} \end{alignat*} 

\subsubsection{Immediate derivative norms}

Without effort we obtain
\begin{alignat*}{5} \normM{(\Aop)(\Aon)\Bo\Rop^{-1}} 
    &\topref{eq:Rop-inverse}{\eq}
    \normM{(\Aop)(\Aon)\Bo\Bo^{-1}(\Aon)^{-1}(\Aop)^{-1}} 
    \\&\eq
    \normM{\id_\Mh}
    \eq
    1  
    \myeqlab{eq:Aop-Aon-Bo-u} \end{alignat*}
\begin{alignat*}{5} \normM{(\Aon)\Bo\Rop^{-1}} 
    &\topref{eq:Rop-inverse}{\eq} 
    \normM{(\Aon)\Bo\Bo^{-1}(\Aon)^{-1}(\Aop)^{-1}} 
    \eq
    \normM{(\Aop)^{-1}} 
    \\&\topref{eq:Aopn-norm}{\lesssim}
    \asy\nv^{-1}  
    \myeqlab{eq:Aon-Bo-u} \end{alignat*} 
\begin{alignat*}{5} \normM{\Bo\Rop^{-1}} 
    &\topref{eq:Rop-inverse}{\eq}
    \normM{\Bo\Bo^{-1}(\Aon)^{-1}(\Aop)^{-1}} 
    \\&\leq 
    \subeq{ \normM{(\Aon)^{-1}} }{ \topref{eq:Aopn-norm}{\lesssim} \asy\nv^{-1} } 
    \subeq{ \normM{(\Aop)^{-1}} }{ \topref{eq:Aopn-norm}{\lesssim} \asy\nv^{-1} } 
    \lesssim
    \asy\nv^{-2}  
    \myeqlab{eq:Bo-u} \end{alignat*} 
\begin{alignat*}{5} \normM{\Rop^{-1}} 
    &\topref{eq:Rop-inverse}{\eq} 
    \normM{\Bo^{-1}(\Aon)^{-1}(\Aop)^{-1}} 
    \\&\leq 
    \subeq{ \normM{\Bo^{-1}} }{ \topref{eq:B-norm}{\lesssim} 1 }
    \subeq{ \normM{(\Aon)^{-1}} }{ \topref{eq:Aopn-norm}{\lesssim} \asy\nv^{-1} } 
    \subeq{ \normM{(\Aop)^{-1}} }{ \topref{eq:Aopn-norm}{\lesssim} \asy\nv^{-1} } 
    \lesssim
    \asy\nv^{-2}  
    \myeqlab{eq:u} \end{alignat*} 

\subsubsection{Other derivative norms}

Linear combinations yield 
\begin{alignat*}{5} 
    \normM{ (\Ao+\mgr)\Bo \Rop^{-1} } 
    &\eq  \normM{ [(\Ao+\mmen)+(\mgr-\mmen)]\Bo \Rop^{-1} } 
    \\&\leq 
    \subeq{ \normM{ (\Ao+\mmen)\Bo \Rop^{-1} } }{ \topref{eq:Aon-Bo-u}{\lesssim} \asy\nv^{-1}  }
    + 
    \subeq{ |\gr-\men| }{ \toprefb{eq:grgsgt}{eq:mepn-bound}{\lesssim} \asy\nv }
    \subeq{ \norm{ \Bo \Rop^{-1} } }{ \topref{eq:Bo-u}{\lesssim} \asy\nv^{-2}  }
    \lesssim \asy\nv^{-1}  
    \myeqlab{eq:Aor-Bo-u}
\end{alignat*} 
\begin{alignat*}{5} \normM{ (\Ao+\mgr)(\Ao+\mgs)\Bo \Rop^{-1} } 
    \eq&
    \normM{ [(\Ao+\mmep)+(\mgr-\mmep)][(\Ao+\mmen)+(\mgs-\mmen)]\Bo \Rop^{-1} }&&
    \\\leq& 
    \subeq{ \normM{ (\Ao+\mmep)(\Ao+\mmen)\Bo \Rop^{-1} } }{ \topref{eq:Aop-Aon-Bo-u}{\lesssim}  1  }
    + 
    \subeq{ |\gr-\mep| }{ \toprefb{eq:grgsgt}{eq:mepn-bound}{\lesssim} \asy\nv }
    \subeq{ \normM{ (\Ao+\mmen)\Bo \Rop^{-1}  } }{ \topref{eq:Aon-Bo-u}{\lesssim} \asy\nv^{-1} }&&
    \\&+ 
    \subeq{ |\gs-\men| }{ \toprefb{eq:grgsgt}{eq:mepn-bound}{\lesssim} \asy\nv }
    \subeq{ \normM{ (\Ao+\mmep)\Bo \Rop^{-1}  } }{ \topref{eq:Aor-Bo-u}{\lesssim} \asy\nv^{-1} }
    + 
    \subeq{ |\gr-\mep| }{ \toprefb{eq:grgsgt}{eq:mepn-bound}{\lesssim} \asy\nv }
    \subeq{ |\gs-\men| }{ \toprefb{eq:grgsgt}{eq:mepn-bound}{\lesssim} \asy\nv }
    \subeq{ \normM{ \Bo \Rop^{-1}  } }{ \topref{eq:Bo-u}{\lesssim} \asy\nv^{-2} }
    &&\lesssim
    1 
    \myeqlab{eq:Aor-Aos-Bo-u}
\end{alignat*} 
To obtain estimates for other derivatives we could try to commute $(\Aopn)^{-1}$ with $\Bo$, or $\Bo^{-1}$ with $\Aopn$.
But this would produce $(\Bo+\moacs1)^{-1}$, which happens to be pathological. However, this can be circumvented:
we have already estimated
\begin{alignat*}{5} 
    (\Ao+\mgr)(\Ao+\mgs)\Bo \Rop^{-1} 
    \topref{eq:arasbt}{\eq}&
    [(\Bo +\moacs2)(\Ao +\mgr-\moacs1)(\Ao +\mgs-\moacs1)
    +\subeq{(-\mgr-\mgs+\moacs3)}{\overset!=0}\Ao
    -\moacs2(\mgr-\moacs1)(\mgs-\moacs1)
    ]\Rop^{-1} 
    \myeqlab{eq:modd}
\end{alignat*} 
Trick: we eliminate the $\Ao$ term, for which we do not have an estimate yet, by specializing \myeqref{eq:modd} 
to $\gr\eq \nv+1$, $\gs\eq 2-\nv$ (for example;
we need $\gr+\gs=3$ and $\gr,\gs\sim\nv$) 
\begin{alignat*}{5}  
    (\Ao+\nvo+\moacs1)(\Ao+\moacs2-\nvo)\Bo \Rop^{-1} 
    &\topref{eq:modd}\eq 
    [ (\Bo+\moacs2)(\Ao+\nvo)(\Ao+\moacs1-\nvo) -\moacs2\nvo(\moacs1-\nvo) ]\Rop^{-1} 
    \myeqlab{eq:ABtrick}
\end{alignat*} 
and now we have an estimate with $\Bo$ on the \emph{left}:
\begin{alignat*}{5} 
    \normM{ (\Bo+\moacs2)(\Ao+\nvo)(\Ao+\moacs1-\nvo)\Rop^{-1}  }
    &\topref{eq:ABtrick}\leq
    \subeq{ \normM{ (\Ao+\nvo+\moacs1)(\Ao+\moacs2-\nvo)\Bo \Rop^{-1}  } }{ \topref{eq:Aor-Aos-Bo-u}{\lesssim}  1  }
    + 
    2|\nv-1||\nv| 
    \subeq{ \norm{\Rop^{-1} } }{ \topref{eq:u}{\lesssim} \asy\nv^{-2} }
    \\&\lesssim
    1 
    \myeqlab{eq:Btwo-An-An}
\end{alignat*} 
Now we can estimate pure-$\Ao$ products:
\begin{alignat*}{5} 
    \normM{ (\Ao+\nvo)(\Ao+\moacs1-\nvo)\Rop^{-1}  }
    &\eq  
    \normM{ (\Bo+\moacs2)^{-1}(\Bo+\moacs2)(\Ao+\nvo)(\Ao+\moacs1-\nvo)\Rop^{-1}  }
    \\&\leq 
    \subeq{ \normM{ (\Bo+\moacs2)^{-1} } }{ \topref{eq:Btwo-norm}{\lesssim} 1 }~
    \subeq{ \normM{ (\Bo+\moacs2)(\Ao+\nvo)(\Ao+\moacs1-\nvo)\Rop^{-1}  } }{ \topref{eq:Btwo-An-An}{\lesssim}  1  }
    \lesssim
    1 
    \myeqlab{eq:An-An}
\end{alignat*} 
and
\begin{alignat*}{5} \normM{ (\Ao+\moacs1-\nvo)\Rop^{-1}  }
    &\eq  
    \normM{ (\Ao+\nvo)^{-1}(\Ao+\nvo)(\Ao+\moacs1-\nvo)\Rop^{-1}  }
    \\&\leq 
    \subeq{\normM{ (\Ao+\nvo)^{-1} }}{\topref{eq:Aoplusn-norm}{\lesssim} \asy\nv^{-1}}~
    \subeq{ \normM{ (\Ao+\nvo)(\Ao+\moacs1-\nvo)\Rop^{-1}  } }{ \topref{eq:An-An}{\lesssim} 1 }
    \lesssim
    \asy\nv^{-1} 
    \myeqlab{eq:An}
\end{alignat*} 
Now linear combinations yield 
\begin{alignat*}{5} \normM{ (\Ao+\mgr) \Rop^{-1}  }
    &\eq 
    \normM{ [(\Ao+\moacs1-\nvo)+(\mgr+\nvo-\moacs1)]\Rop^{-1}  }
    \\&\leq
    \subeq{ \normM{ (\Ao+\moacs1-\nvo) \Rop^{-1}  } }{ \topref{eq:An}{\lesssim} \asy\nv^{-1} }
    + 
    \subeq{ |\mgr+\nvo-\moacs1| }{ \topref{eq:grgsgt}{\lesssim}\asy\nv }
    \subeq{ \normM{ \Rop^{-1}  } }{ \topref{eq:u}{\lesssim} \asy\nv^{-2} }
    \lesssim
    \asy\nv^{-1} 
    \myeqlab{eq:Ar}
\end{alignat*} 
and
\begin{alignat*}{5} \normM{ \Ao\Ao \Rop^{-1}  }
    \eq&
    \normM{ [(\Ao+\nvo)-\nvo][(\Ao+\moacs1-\nvo)+(\nvo-\moacs1)]\Rop^{-1}  } &&
    \\\leq&
    \subeq{ \normM{ (\Ao+\nvo)(\Ao+\moacs1-\nvo)\Rop^{-1}  } }{\topref{eq:An-An}{\lesssim} 1 }
    + |\nv| \subeq{ \normM{ (\Ao+\moacs1-\nvo)\Rop^{-1}  } }{\topref{eq:Ar}{\lesssim}\asy\nv^{-1} } &&
    \\&+ |\nv-1| \subeq{ \normM{ (\Ao+\nvo)\Rop^{-1}  } }{\topref{eq:Ar}{\lesssim}\asy\nv^{-1} }
    + |\nv| |\nv-1| \subeq{ \normM{ \Rop^{-1}  } }{\topref{eq:u}{\lesssim}\asy\nv^{-2} } &&
    \lesssim
    1 
    \myeqlab{eq:As-Ar}
\end{alignat*} 
and the rest is straightforward:
\begin{alignat*}{5} \normM{\Bo\Ao \Rop^{-1} } 
    &\topref{eq:BAAB}{\eq} 
    \normM{(\Ao\Bo+\Bo-\Ao)\Rop^{-1} } 
    \leq
    \subeq{ \normM{\Ao\Bo \Rop^{-1} } }{\topref{eq:Aor-Bo-u}{\lesssim} \asy\nv^{-1} }
    +\subeq{ \normM{\Bo \Rop^{-1} } }{\topref{eq:Bo-u}{\lesssim} \asy\nv^{-2} }
    +\subeq{ \normM{\Ao \Rop^{-1} } }{\topref{eq:Ar}{\lesssim} \asy\nv^{-1} }
    \lesssim
    \asy\nv^{-1} 
    \myeqlab{eq:Bo-Ao-u}
\end{alignat*} 
\begin{alignat*}{5} 
    \normM{\Ao\Bo\Ao \Rop^{-1} }
    &\topref{eq:BAAB}{\eq}
    \normM{\Ao(\Ao\Bo+\Bo-\Ao)\Rop^{-1} }
    \leq
    \subeq{\normM{\Ao\Ao\Bo \Rop^{-1} }}{\topref{eq:Aor-Aos-Bo-u}{\lesssim} 1 }
    +\subeq{\normM{\Ao\Bo \Rop^{-1} }}{\topref{eq:Aor-Bo-u}{\lesssim}\asy\nv^{-1} }
    +\subeq{\normM{\Ao\Ao \Rop^{-1} }}{\topref{eq:As-Ar}{\lesssim} 1 }
    \lesssim
    1 
    \myeqlab{eq:ABA-u}
\end{alignat*} 
\begin{alignat*}{5} 
    \normM{\Bo\Ao\Ao \Rop^{-1} }
    &\topref{eq:BAAB}{\eq}
    \normM{(\Ao\Bo+\Bo-\Ao)\Ao \Rop^{-1} }
    \leq
    \subeq{ \normM{\Ao\Bo\Ao \Rop^{-1} } }{ \topref{eq:ABA-u}{\lesssim}  1  }
    +
    \subeq{ \norm{\Bo\Ao \Rop^{-1} } }{ \topref{eq:Bo-Ao-u}{\lesssim} \asy\nv^{-1}  }
    +
    \subeq{ \normM{\Ao\Ao \Rop^{-1} } }{ \topref{eq:As-Ar}{\lesssim}  1  }
    \lesssim
    1
    \myeqlab{eq:BAA-u}
\end{alignat*} 

\subsection{$\Lopn$ inversion}

\mylabel{section:Lop-inversion}%
$\Rop$ is an isomorphism on $\Rop^{-1}[\Mh]$ (induced norm) into $\Mh$. 
We have to show that 
\begin{alignat*}{5} \Eop \topref{eq:lstff-new}{\eq} (\moacs2\moacs\mv-\moacs1)(\Ao-\Bo):\Rop^{-1}[\Mh]\rightarrow\Mh  \end{alignat*} 
is small in that same norm, which is equivalent to showing that $\Eop\Rop^{-1}:\Mh\rightarrow\Mh$ is small.

For $\nv=0$ we observe that $\Ao-\Bo\toprefb{eq:mdU-fou}{eq:mdta-fou-a}{=}(\etao-\nvo)\peta-\etao\peta=0$ so that $\Eop=0$. 

Consider $\nv\neq 0$: 
\begin{alignat*}{5} 
    \normM{ \Eop\Rop^{-1} } 
    &\lesssim 
    \subeq{ \normM{\Ao \Rop^{-1} } }{ \topref{eq:Ar}{\lesssim} \asy\nv^{-1} } 
    +
    \subeq{ \normM{\Bo\Rop^{-1}} }{ \topref{eq:Bo-u}{\lesssim} \asy\nv^{-2} }
    \lesssim
    \asy\nv^{-1} 
\end{alignat*} 
Therefore $\normM{\Eop\Rop^{-1}}<1$ if we pick $\Nv$ sufficiently large (since $\nv\in\Nv\Z\setdiff\{0\}$ means $|\nv|\geq\Nv$). 
Then $\Rop-\Eop$ is also an isomorphism, with inverse
\begin{alignat*}{5} 
    \Lopn^{-1} = (\Rop-\Eop)^{-1} 
    = \Rop^{-1} \subeq{ \sum_{k=0}^\infty (\Eop\Rop^{-1})^k }{=\Sop} 
    \myeqlab{eq:Lopex}
\end{alignat*} 
\begin{alignat*}{5} 
    \normM{\Sop} \leq \sum_{k=0}^\infty\normM{(\Eop\Rop^{-1})^k} \leq \frac{1}{1-\normM{\Eop\Rop^{-1}}} < \infty \myeqlab{eq:ERnorm}
\end{alignat*} 
so the estimates for $\Rop^{-1}f$ from the previous section yield 
\begin{alignat*}{5} 
    \normM{\Lopn^{-1}} 
    &\topref{eq:Lopex}{\eq}
    \normM{\Rop^{-1}\Sop} 
    \leq 
    \subeq{ \normM{\Rop^{-1}} }{ \topref{eq:u}{\lesssim} \asy\nv^{-2} }
    \subeq{ \normM{\Sop} }{ \topref{eq:ERnorm}{\lesssim} 1 }
    \lesssim
    \asy\nv^{-2}
    \myeqlab{eq:I}
    \\
    \normM{\Bo \Lopn^{-1} } 
    &\topref{eq:Lopex}{\eq}
    \normM{ \Bo \Rop^{-1} \Sop }
    \leq
    \subeq{ \normM{ \Bo \Rop^{-1} } }{ \topref{eq:Bo-u}{\lesssim} \asy\nv^{-2} }
    \subeq{ \normM{ \Sop } }{ \topref{eq:ERnorm}{\lesssim} 1 }
    \lesssim
    \asy\nv^{-2}
    \myeqlab{eq:B}
\end{alignat*} 
and similarly
\begin{alignat*}{5} 
    \normm{\Ao \Lopn^{-1} },\norm{\Ao\Bo \Lopn^{-1} },\norm{\Bo\Ao \Lopn^{-1} } &\toprefc{eq:Ar}{eq:Aor-Bo-u}{eq:Bo-Ao-u}{\lesssim} \asy\nv^{-1} \quad,
    \myeqlab{eq:A-AB-BA}
    \\
    \normm{\Ao\Ao \Lopn^{-1} },\norm{\Ao\Ao\Bo \Lopn^{-1} },\norm{\Ao\Bo\Ao \Lopn^{-1} },\norm{\Bo\Ao\Ao \Lopn^{-1} } 
    &\toprefd{eq:As-Ar}{eq:Aor-Aos-Bo-u}{eq:ABA-u}{eq:BAA-u}{\lesssim}  1 \quad.
    \myeqlab{eq:AA-AAB-ABA-BAA}
\end{alignat*}

\renewcommand{\normm}[1]{\norma{#1}{\Mh}}
\renewcommand{\normM}[1]{\norma{#1}{[\Mh]}}

\section{Existence and local uniqueness and stability}

\mylabel{section:existence}

\subsection{Function spaces}

The values of $\Fmap$ (see \myeqref{eq:curl-scaled-new}) will be in the \defm{defect space}
\begin{alignat*}{5} \fdom := \big\{ \fff:\R\times\T\rightarrow\R ~\big|~ \ffff \in \loner(\Mh) \big\} \myeqlab{eq:fdom} \quad; \end{alignat*} 
we also consider its complexification
\begin{alignat*}{5} \fdomC := \CC\tensor\fdom = \big\{ \fff:\R\times\T\rightarrow\CC ~\big|~ \ffff \in \loner(\Mh) \big\} = \loner(\Mh)^\vee \quad, \myeqlab{eq:fdomC} \end{alignat*} 
both with norm 
\begin{alignat*}{5} \norma{\fff}{\fdom} := \norma{\ffff}{\loner(\Mh)} = \sum_{\nv\in\Z}\norma{\ffff(\cdot,\nv)}{\Mh}\quad. \end{alignat*} 
$\Mh$ is a closed subalgebra of $\Ms$ and $\loner(X)$ a closed subalgebra of $\lone(X)$, both with convolution $*$ as
product, so
$\fdom$ is a Banach algebra with pointwise multiplication $\cdot$ as product. 
Moreover, by Riemann-Lebesgue 
\begin{alignat*}{5} \fdom\embed\Cubk0(\R\times\T;\R) \myeqlab{eq:fdom-Cub} \end{alignat*} 
The elements of $\fdom$ are $\frac{2\pi}{\Nv}$-periodic in $\uan$.

Define the subspace
\begin{alignat*}{5} \Vortdom:=\{\fff\in\fdom~|~\text{$\fff(\ta,\uan)$ constant in $\ta$} \}\quad; \myeqlab{eq:Vorfdom} \end{alignat*}
we seek
\begin{alignat*}{5} \Vort \in \Vortdom &\embed \Cubk0(\T) \myeqlab{eq:Vortdom-Cub} \end{alignat*} 

Abbreviate $\Lop:=\frac{\partial\Fmap}{\partial\sstf}(\sstfz,\sVortz)^{-1}$ as in \myeqref{eq:linearization}.
Define 
\begin{alignat*}{5} \stfdomC:=\Lop^{-1}[\fdomC]  \myeqlab{eq:stfdomC} \end{alignat*}
We have shown from \myeqref{eq:linearization} to \myeqref{eq:lstff-new} 
that 
\begin{alignat*}{5} (\Lop\lsstf)^\wedge(\cdot,\nv)=\Lopn(\lsstf^\wedge(\cdot,\nv))\quad. \end{alignat*} 
Thus
\begin{alignat*}{5} \stfdomC 
&= \big\{ \sstf~|~\Lop\sstf\in\fdomC \big\} 
\topref{eq:fdomC}{=} 
\big\{ \sstf~|~\sstf^\wedge\in\Lopn^{-1}[\loner(\Mh)] \big\} 
\end{alignat*} 
We have obtained in Section \myref{section:inversion} that $\Lopn:\Lopn^{-1}[\Mh]\rightarrow\Mh$ is bijective 
for every $\nv\in\Nv\Z$ (by taking $\Nv$ sufficiently large), 
so 
$\Lop:\stfdomC\rightarrow\fdomC$ is also bijective. If we endow $\stfdomC=\Lop^{-1}[\fdomC]$ with the induced norm, then $\Lop$ is isometry.

Define 
\begin{alignat*}{5} \stfdom:=\{\sstf:\R\times\T\rightarrow\utilde{~\R~}~|~\sstf\in\stfdomC\} \quad\myeqlab{eq:stfdom}  \end{alignat*} 
$\Lop$ maps real functions into real functions (see \myeqref{eq:linearization}), so $\Lop$ restricted to $\stfdom\rightarrow\fdom$ is still an isometry.
($\stfdom,\fdom$ are closed subspaces of $\stfdomC,\fdomC$, hence also Banach.)

\subsection{Small neighbourhoods}

The estimates \myeqref{eq:I}, \myeqref{eq:B}, \myeqref{eq:A-AB-BA} and \myeqref{eq:AA-AAB-ABA-BAA} show 
\begin{alignat*}{5} \forall k=0..,2,j=0..,k,i=0,1:~\text{$\Ao^j\nvo^{k-j}\Bo^i$ is continuous on $\Lopn^{-1}[\Mh]$ into $\Mh$.} \end{alignat*} 
Inverse transform (see \myeqref{eq:mdU-fou}, \myeqref{eq:mdta-fou}, \myeqref{eq:mdu-fou}) yields, 
by definition of $\fdom,\stfdom$, that
\begin{alignat*}{5} \forall k=0...2,j=0...k,i=0,1:~\text{$\mdU^j\mdu^{k-j}\mdta^i$ is continuous on $\stfdom$ into $\fdom$,} \myeqlab{eq:regu} \end{alignat*} 
hence by \myeqref{eq:fdom-Cub}
\begin{alignat*}{5} \forall k=0...2,j=0...k,i=0,1:~\text{$\mdU^j\mdu^{k-j}\mdta^i$ is continuous on $\stfdom$ into $\Cubk0(\R\times\T)$} \myeqlab{eq:reguC} \end{alignat*} 
Therefore 
\begin{alignat*}{5} 
    \srad &\topref{eq:crad-srad}{=} (\frac{\mdta\sstf}{-\mv})^{-\frac12} \in \Cubk0(\R\times\T) \myeqlab{eq:srad-cub} \\
    \svorf &\topref{eq:cvorf-svorf}{=} (\mdU\sstf)^{-\frac1{2\mv}} \in \Cubk0(\R\times\T) \myeqlab{eq:svorf-cub} 
\end{alignat*} 
For the rest of this paper,
let $a\lesssim b$ denote $|a|\leq Cb$ for some constant $C<\infty$ independent of $\sstf$ and $\ta,\uan$ (or equivalent coordinates), 
and let $a\sim b$ imply $b\lesssim a$ as well. We write $\lesssim^\tv$ to indicate $C$ is also independent of $\tv$.

$\sVortz\topref{eq:Vortz}{=}\frac{2\mv-1}{\mv}>0$, so we can pick $\frad$ so small that for $\Vort\in\Vortnbh$ (ball of radius $\frad>0$ in $\Vortdom$ norm), 
\myeqref{eq:Vortdom-Cub} implies
\begin{alignat*}{5} \Vort &\sim 1 \myeqlab{eq:nonzero-Vort} \end{alignat*} 
Likewise, 
$-(\mdU+\moacs1)\mdta\sstfz\topref{eq:vorfbase}{=}2\mv>0$ and $\mdU\sstfz\topref{eq:mdU-sstfz}{=}1>0$ and 
$-\mdta\sstfz\topref{eq:mdta-sstfz}{=}1>0$, so we can pick $\brad$ so small that for $\sstf\in\stfnbh$, \myeqref{eq:fdom-Cub} implies 
\begin{alignat*}{5} 
    -(\mdU+\moacs1)\mdta\sstf &\sim 1 \myeqlab{eq:nonzero-mdUonemdta}\\
    \mdU\sstf &\sim 1 \myeqlab{eq:nonzero-mdU}\\
    -\mdta\sstf &\sim 1 \myeqlab{eq:nonzero-mdta}
\end{alignat*}
and therefore 
\begin{alignat*}{5}  
    \srad &\topref{eq:crad-srad}{=} (\frac{\mdta\sstf}{-\mv})^{\frac12} \topref{eq:nonzero-mdta}{\sim} 1 \myeqlab{eq:nonzero-srad} \\
    \svorf &\topref{eq:cvorf-svorf}{=} (\mdU\sstf)^{-\frac1{2\mv}} \topref{eq:nonzero-mdU}{\sim} 1 \myeqlab{eq:nonzero-svorf} \\
    \svort &\topref{eq:cvort-svort}{=} \svorf\Vort \toprefb{eq:nonzero-svorf}{eq:nonzero-Vort}{\sim} 1 \myeqlab{eq:nonzero-svort} \\
\end{alignat*} 
In any case, \myeqref{eq:fdom-Cub} yields  
\begin{alignat*}{5} 
    \mdT\sstf  &\toprefb{eq:mdT}{eq:reguC}{\lesssim} 1 \myeqlab{eq:mdT-sstf-bounded} \\
    \mda\sstf  &\toprefb{eq:mda}{eq:reguC}{\lesssim} 1 \myeqlab{eq:mda-sstf-bounded} 
\end{alignat*} 
and (using the power rule $\mdU(f^s)=\ta\pU(f^s)=sf^{s-1}\ta\pU f=sf^{s-1}\mdU f$) 
\begin{alignat*}{5} \mdU\svorf &\toprefb{eq:cvorf-svorf}{eq:reguC}{\lesssim} 1 \myeqlab{eq:mdU-svorf-bounded} \end{alignat*}     

\subsection{Existence}

Consider $\Fmap$ in \myeqref{eq:curl-scaled-new}. By tracing back its constituents it can be verified that $\Fmap$ is a function of $\sVort$ and of those derivatives of $\sstf$ that appear in \myeqref{eq:regu}. 
$\Fmap$ is composed from all these by linear combination, pointwise multiplication,
quotients with denominator $(\mdU+\moacs1)\mdta\sstf$, and fractional powers of $\mdU\sstf$ (see above). 
The last two can be expanded at $\sstfz$ into 
power series that converge absolutely on $\stfnbh$ for $\brad$ sufficiently small, due to \myeqref{eq:nonzero-mdUonemdta} and \myeqref{eq:nonzero-mdta}. 
Therefore
\begin{alignat*}{5} \Fmap\in\Cone(\stfnbh\times\Vortnbh;\fdom) \myeqlab{eq:Fmap-Cone} \end{alignat*} 
In addition 
$\partial\Fmap/\partial\sstf(\sstfz,\Vortz):\stfdom\rightarrow\fdom$ was shown to be a linear isomorphism,
and $\Fmap(\sstfz,\Vortz)=\zeroacc\Fmap\topref{eq:background-solves}{=}0$.
Thus we may invoke the implicit function theorem for Banach spaces. 
It shows that for $\frad$ sufficiently small
there is a $\Cone$ map $H:\Vortnbh\rightarrow\stfnbh$ so that $\Fmap(\sstf,\Vort)\eq 0$ for $\sstf\eq H(\Vort)$ for all $\Vort\in\Vortnbh$.
\begin{alignat*}{5} \text{For the remainder of the paper we set $\sstf\defeq H(\sVort)$.} \end{alignat*}

\subsection{Reversing coordinate changes}

\mylabel{section:coordchange}

Note
\begin{alignat*}{5} \pta\topref{eq:pU-def}{=}\pu-\pU \myeqlab{eq:pta-pu-pU} \end{alignat*}
so after taking an initial $\pta$ (or none) of $\cstf$, \myeqref{eq:regu} allows taking another $\pu^2$ or $\pu\pta$ or $\pta^2$, 
with result continuous at any\footnote{not $\ta\leq 0$ due to omitting the $\tao$ from $\mdU,\mdta$} $(\ta,\uan)\in\boi0\infty\times\T$.
Hence 
\begin{alignat*}{5} \crad,\lrv \toprefb{eq:crad-srad}{eq:rad-lrv}{\in} \Ck2(\Rplus\times\T) \myeqlab{eq:rad-regu} \end{alignat*} 
(after the $\pta$ in \myeqref{eq:crad-srad} we can take two more $\pta$ or $\pu$).

\begin{alignat*}{5} 
    \crad &\topref{eq:crad-srad}{=} \ta^{-\mv}\srad > 0 \quad, \myeqlab{eq:nonzero-crad} 
\end{alignat*} 
so $\lrv=\log\crad$ is well-defined, and
\begin{alignat*}{5} 
    \lrv_\Uan
    \topref{eq:pU-lrv-scaled}{=} 
    \ta^{-1}\frac{(\mdU+\moacs1)\mdta\sstf}{\mdta\sstf}
    \toprefb{eq:nonzero-mdUonemdta}{eq:nonzero-mdta}{\sim}
    \ta^{-1}
    >0
    \myeqlab{eq:aU-good}
\end{alignat*} 
We consider coordinates $(\polb,\Uan)$ (see Figure \myref{fig:polb-Uan})
\begin{alignat*}{5} \polb = \ta+\uan \csep \Uan = \uan \myeqlab{eq:polbUan} \end{alignat*} 
so that $\pU$ is a \emph{radial} derivative and $\polb$ angular. $\ta>0$ corresponds to $\Uan<\polb$.

Fix $\polb$. Integrating 
\begin{alignat*}{5} \lrv_\Uan \sim
    \ta^{-1} 
    =
    (\polb-\Uan)^{-1} > 0 \end{alignat*} 
yields that $\Uan\mapsto\lrv(\polb,\Uan)$ is strictly increasing, and 
\begin{alignat*}{5} 
    \lrv(\polb,\Uan) \overset{\Uan\nearrow\polb}{\rightarrow} \infty 
    \csep
    \lrv(\polb,\Uan) \overset{\Uan\searrow-\infty}{\rightarrow} -\infty
\end{alignat*} 
so $\boi{-\infty}{\polb} \ni \Uan \mapsto \lrv(\polb,\Uan) \in \R$ is bijective.

Moreover $\pola\topref{eq:pola}{=}\ta+\uan=\polb$. Therefore the map 
\begin{alignat*}{5} \bigcup_{\polb\in\R}(\{\polb\}\times\boi{-\infty}{\polb})\ni(\polb,\Uan)\mapsto(\ta,\uan)\mapsto(\lrv,\pola)\in\boi0\infty\times\R  \end{alignat*} 
bijective,
hence a $\Ck2$ diffeomorphism by regularity estimates \myeqref{eq:rad-regu}. 
The transform from $(\lrv,\pola)$ to $\cvx$ is also a diffeomorphism (modulo periodicity).

\section{Euler solution}

\subsection{Classical solution}

\begin{alignat*}{5} 
    \stf 
    \topref{eq:stf-cstf}{=} 
    \tv^{2\mv-1}\cstf
    \topref{eq:cstf-sstf}{=}
    (\frac{\tv}{\ta})^{2\mv-1} \sstf
    \topref{eq:regu}{\in} \Ck2(\Rplus\times\T)
\end{alignat*} 
Since the coordinate transforms are $\Ck2$, we get
\begin{alignat*}{5} \stf=\stf(\vx,\tv) \in \Ck2((\R^2\backslash\{0\})\times\boi0\infty), \myeqlab{eq:stf-regu} \\
    \impl\quad \vv^{(\vx)} \topref{eq:stf}{\defeq} \nabla_\vx^\perp\stf 
    \in \Ck1((\R^2\backslash\{0\})\times\boi0\infty) \myeqlab{eq:vv-regu} \end{alignat*} 
We have a classical solution of
\begin{alignat*}{5} \nabla\cdot\vv = 0 \quad, \myeqlab{eq:divvv} \end{alignat*} 
and 
by tracing the steps from \myeqref{eq:curl-scaled-new} back to the original curl equation \myeqref{eq:vort-cvort},
we have obtained a classical solution of 
\begin{alignat*}{5} \Delta\stf = \vort \myeqlab{eq:vorttt} \end{alignat*} 
as well. Hence 
\begin{alignat*}{5} \vort \toprefb{eq:vorttt}{eq:stf-regu}{\in} \Ck0((\R^2\backslash\{0\})\times\boi0\infty) \end{alignat*}
and 
\begin{alignat*}{5} \nabla\times\vv \toprefb{eq:vorttt}{eq:vv-regu}{=} \vort\quad. \end{alignat*}  
Tracing back from \myeqref{eq:vort-stf-ta} to the original vorticity equation \myeqref{eq:vort-div}, 
and observing that the distributional definition of the outermost derivatives has been respected, 
we obtain a weak (continuous) solution $\vort$ of
\begin{alignat*}{5} 
    0 
    = \pt\vort + \nabla\cdot(\vort\vv) 
    \toprefb{eq:divvv}{eq:vorttt}{=} 
    \nabla\times(\pt\vv+\nabla\cdot(\vv\otimes\vv))
    \myeqlab{eq:curl-veq}
\end{alignat*} 
Combining these three equations shows that in some neighbourhood of each $(\vx,t)\in(\R^2\backslash\{0\})\times\Rplus$ 
there is a $\Ck1$ (in particular) scalar $\prs$ so that
\begin{alignat*}{5} 
    0 &= \pt\vv+\nabla\cdot(\vv\otimes\vv)+\nabla\prs 
    \myeqlab{eq:classical-locally}
\end{alignat*} 

\subsection{Convergence to initial data}

\begin{alignat*}{5} \rad &\topref{eq:vx-cvx}{=} \tv^{\mv}\crad \topref{eq:nonzero-crad}{=} (\frac{\tv}{\ta})^\mv\srad 
    \topref{eq:nonzero-srad}{\sim^\tv}
    (\frac{\tv}{\ta})^{\mv}
    > 
    0\quad, \myeqlab{eq:nonzero-rad} \\
    \vort 
    &\topref{eq:vort-cvort}{=} 
    \tv^{-1}\cvort 
    \topref{eq:cvort-svort}{=} 
    (\frac{\tv}{\ta})^{-1}\svort 
    \topref{eq:nonzero-svort}{\sim^\tv}
    (\frac{\tv}{\ta})^{-1}
    \topref{eq:nonzero-rad}{\sim^\tv}
    \rad^{-\frac1\mv}
    \topref{eq:nonzero-svort}{>} 0
    \quad, \myeqlab{eq:nonzero-vort} 
\end{alignat*} 
Consider a fixed $\vx\neq 0$ and therefore fixed $\pola\topref{eq:crad-pola-def}{=}\measuredangle\cvx=\measuredangle\vx$. 
Then $t\searrow 0$ yields 
\begin{alignat*}{5} 
    \ta &\topref{eq:nonzero-rad}{\sim^\tv} \tv\rad^{-\frac1\mv} \rightarrow 0 \quad\text{and} \quad \uan \topref{eq:pola}{=} \pola-\ta\rightarrow\pola 
    \myeqlab{eq:tbconv}
\end{alignat*} 
Abbreviate
\begin{alignat*}{5} \vrd \defeq \vort \rad^{\frac1\mv} 
    \toprefb{eq:nonzero-vort}{eq:nonzero-rad}{\eq}
    (\frac{\tv}{\ta})^{-1} \svort 
    \cdot
    \big((\frac{\tv}{\ta})^{\mv} \srad\big)^{\frac1\mv} 
    =
    \svort\srad^{\frac1\mv}
    \topref{eq:cvort-svort}{=}
    \sVort\svorf\srad^{\frac1\mv}
    \myeqlab{eq:vrd}
\end{alignat*} 
Since $\Vort,\svorf,\srad$ are continuous\footnote{They are constructed as functions of $\ta\in\R$, not just $\ta\in\boi0\infty$, and $\uan\in\T$.} at $\ta=0$ (by \myeqref{eq:Vortdom-Cub}, \myeqref{eq:srad-cub}, \myeqref{eq:svorf-cub}), so is $\vrd$, hence
\begin{alignat*}{5} \vort 
    \toprefb{eq:vrd}{eq:tbconv}{\rightarrow}
    \rad^{-\frac1\mv} \subeq{\vrd(\ta=0,\uan=\pola)}{=:\vrdini(\pola)}
    \quad\text{as $\tv\searrow 0$}
    \myeqlab{eq:vrd-lim-infinity} \end{alignat*} 
with locally uniform convergence; this is \myeqref{eq:vort-ini-cond}. 
(We analyze in section \myref{section:ini-data} how $\vrdini$ depends on $\sVort$.)

Then using $\vv^{(\vx)}=\nabla_\vx^\perp\stf$ with $\Delta\stf=\vort$  we obtain 
\begin{alignat*}{5} \text{$\vv^{(\vx)}$ converges locally uniformly in any $\vx\neq 0$ as $t\searrow 0$.} \myeqlab{eq:vv-convergence} \end{alignat*}

\subsection{Weak solution on entire domain}

\newcommand{\tha}[1]{\chi_{#1}}
\newcommand{\tho}{\tha1}
\newcommand{\thparm}{\delta}
\newcommand{\thd}{\tha\thparm}
\newcommand{\wpot}{\alpha}
\newcommand{\wpotd}{\wpot^\thparm}
\newcommand{\qvd}{{\qv'}}
\newcommand{\vwd}{\vw^\thparm}

To check whether $\vv^{(\vx)}$ is a weak solution at $\vx=0$ as well, we first analyze its asymptotics:
\begin{alignat*}{5} 
    \vv^{(\vx)}
    &\topref{eq:vv-cvv}{\eq}
    t^{\mv-1}\cvv^{(\cvx)}
    \topref{eq:cvv-cstf}{\eq}
    t^{\mv-1}\skw\nabla_\cvx\cstf
    \eq
    t^{\mv-1}\skw\nabla_\cvx\Tav^T\nabla_\Tav\cstf
    \topref{eq:DDxa}{\eq}
    t^{\mv-1}
    \skw
    \crad^{-1}
    \begin{bmatrix}
        \cos\pola & -\sin\pola \\
        \sin\pola & \cos\pola
    \end{bmatrix}
    \begin{bmatrix}
        \cstf_\lrv \\
        \cstf_\pola
    \end{bmatrix}
    \\&\toprefc{eq:crad-srad}{eq:mdT}{eq:mda}{\eq}
    t^{\mv-1} \skw (\ta^{-\mv}\srad)^{-1}
    \begin{bmatrix}
        \cos\pola  & -\sin\pola \\
        \sin\pola & \cos\pola
    \end{bmatrix}
    \ta^{1-2\mv}
    \begin{bmatrix}
        \mdT \sstf \\
        \mda \sstf
    \end{bmatrix}
    \\&\eq
    (\subeq{\frac{t}{\ta}}{\topref{eq:nonzero-rad}{\sim}\rad^{\frac1\mv}})^{\mv-1} 
    \skw 
    \subeq{\srad^{-1}}{\topref{eq:nonzero-srad}{\sim}1}
    \begin{bmatrix}
        \cos\pola & -\sin\pola \\
        \sin\pola & \cos\pola
    \end{bmatrix}
    \subeq{
        \begin{bmatrix}
            \mdT \sstf \\
            \mda \sstf
        \end{bmatrix}
    }{\toprefb{eq:mdT-sstf-bounded}{eq:mda-sstf-bounded}{\lesssim 1}}
    \lesssim
    \rad^{1-\frac1\mv}
    \quad\text{uniformly in $\tv>0$.}
    \myeqlab{eq:vv-svv-rad}
\end{alignat*}    
We will need a $\qv>2$ so that
\begin{alignat*}{5} |\vv^{(\vx)}(\cdot,\tv)|^2 \lesssim \rad^{2(1-\frac1\mv)} \in \Leba{\qv}_\loc(\R^2) \quad\text{uniformly in $\tv\geq 0$.} \myeqlab{eq:vvq-Lqloc} \end{alignat*} 
It exists if and only if $2\cdot 2(1-\frac1\mv)>-2$ which is equivalent to $\mv>\frac23$, and then
\begin{alignat*}{5} |\vv^{(\vx)}(\cdot,\tv)| \in \Leba{2\qv}_\loc(\R^2) \embed \Leba{\qv}_\loc(\R^2) \quad\text{uniformly in $\tv\geq 0$.} \myeqlab{eq:vv-Lqloc} \end{alignat*}

Let $\vw\in\Cinf_c(\R^2\times\roi0\infty)$ be a divergence-free test function. 
Then $\vw=\nabla^\perp\wpot$ for some scalar $\wpot\in\Cinf(\R^2\times\roi0\infty)$. 
We may assume that $\wpot(0,\tv)=0$ for all $\tv\geq 0$, so that
\begin{alignat*}{5} \wpot = O(\rad) \csep \nabla\wpot = O(1) \csep \nabla^2\wpot = O(1) \quad\text{as $\rad\searrow 0$.} \myeqlab{eq:phi-asy} \end{alignat*} 
To use that we already know $\vv^{(\vx)}$ is a classical solution in $\vx\neq 0$, we cut the origin out of the support of $\wpot$:
pick a radial $\tho\in\Cinf(\R^2)$ with $\tho=1$ on $\ballrc10$ and $\tho=0$ on $\R^2\setdiff\ballrc20$.
Define $\thd(\vx):=\tho(\vx/\delta)$, then $\thd=1$ on $\ballrc\delta0$ and 
\begin{alignat*}{5} 
    \supp\thd &\subset \cballrc{2\delta}0 \myeqlab{eq:thd-supp} \\
    \supp(1-\thd) &\subset \R^2\setdiff\ballrc{\delta}0 \myeqlab{eq:one-thd-supp} \\
    \supp\nabla\thd &\subset \cballrc{2\delta}0\setdiff\ballrc{\delta}0 \myeqlab{eq:supp-Dthd}
    \\
    \thd &= O(1) \csep \nabla\thd = O(\delta^{-1}) \csep \nabla^2\thd = O(\delta^{-2}) \quad\text{as $\delta\searrow 0$,} \end{alignat*} 
so \myeqref{eq:phi-asy} combined with \myeqref{eq:thd-supp} yields
\begin{alignat*}{5} \wpot = O(\delta) \csep \nabla\wpot=O(1) \csep \nabla^2\wpot=O(1) \quad\text{on $\supp\thd$} \end{alignat*} 
and therefore
\begin{alignat*}{5} 
    \wpot\nabla\thd = O(\delta)O(\delta^{-1})=O(1) 
    \csep
    \thd\nabla\wpot &= O(1)O(1) = O(1) \quad,
    \\
    \wpot\nabla^2\thd = O(\delta)O(\delta^{-2}) = O(\delta^{-1})
    \csep
    \nabla\wpot\nabt\thd &= O(1)O(\delta^{-1}) = O(\delta^{-1})\quad,
    \\
    \thd \nabla^2\wpot &= O(1)O(1) = O(1)
    \quad\text{as $\delta\searrow0$,}
\end{alignat*} 
so
\begin{alignat*}{5} \nabla(\thd\wpot) = O(1) \csep \nabla^2(\thd\wpot) = O(\delta^{-1}) \quad\text{as $\delta\searrow0$.} \end{alignat*} 
That means
\begin{alignat*}{5} \int_{\R^2} |\nabla^2(\thd\wpot)|^\qvd d\xv \leq \int_{\ballrc{2\delta}0} O(\delta^{-\qvd}) d\xv = O(\delta^{2-\qvd}) \end{alignat*} 
which converges to $0$ iff $\qvd<2$. Hence
\begin{alignat*}{5} 
    \nabla^\perp((1-\thd)\wpot) &\rightarrow \nabla^\perp\wpot = \vw \quad\text{and}
    \\
    \nabla^\perp\nabt((1-\thd)\wpot) &\rightarrow \nabla^\perp\nabt\wpot = \nabt\vw
    \quad\text{in $\Leba{\qvd}(\R^2)$ as $\delta\searrow 0$, uniformly in $\tv>0$} 
    \myeqlab{eq:phi-weakLone}
\end{alignat*} 
The same results hold with $\wpot$ replaced by $\pt\wpot$ (note $\wpot(0,\tv)=0$ implies $\pt\wpot(0,\tv)=0$).

This is enough to obtain a solution at $\vx=0$ as well (write $\vv=\vv^{(\vx)}$ now): for $\tau>0$,
\begin{alignat*}{5} 
    0
    \toprefb{eq:curl-veq}{eq:one-thd-supp}{=}&
    \int_\tau^\infty 
    \int_{\R^2} 
    \nabla\times(\pt\vv+\nabla\cdot(\vv\otimes\vv)) 
    ~
    (1-\thd)\wpot
    d\vx
    ~dt
    \\=&
    -
    \int_\tau^\infty 
    \int_{\R^2} 
    (\pt\vv+\nabla\cdot(\vv\otimes\vv)) 
    \cdot
    \nabla^\perp((1-\thd)\wpot)
    d\vx
    ~dt
    \\=&
    \Big(
    \int_{\R^2} 
    \vv
    \cdot
    \nabla^\perp((1-\thd)\wpot)
    d\vx
    \Big)_{|t=\tau}
    \\&+
    \int_\tau^\infty 
    \int_{\R^2} 
    \vv
    \cdot
    \nabla^\perp((1-\thd)\pt\wpot)
    +
    \vv^2:
    \nabla^\perp\nabt((1-\thd)\wpot)
    d\vx
    ~dt
    \\\topref{eq:vv-convergence}{\underset{\tau\searrow 0}{\rightarrow}}&
    \Big(
    \int_{\R^2} 
    \subeq{\vv}{\topref{eq:vv-Lqloc}{\in}\Leba{\qv}_\loc}
    \cdot
    \subeq{
        \nabla^\perp((1-\thd)\wpot)
    }{\topref{eq:phi-weakLone}{\underset{\Leba{\qvd}}{\rightarrow}}\vw}
    d\vx
    \Big)_{|t=0}
    \\&+
    \int_{0}^\infty 
    \int_{\R^2} 
    \subeq{\vv}{\topref{eq:vv-Lqloc}{\in}\Leba{\qv}_\loc}
    \cdot
    \subeq{
        \nabla^\perp((1-\thd)\pt\wpot)
    }{\topref{eq:phi-weakLone}{\underset{\Leba{\qvd}}{\rightarrow}}\pt\vw}
    +
    \subeq{\vv^2}{\topref{eq:vv-Lqloc}{\in}\Leba{\qv}_\loc}
    :
    \subeq{\nabla^\perp\nabt((1-\thd)\wpot)}{\topref{eq:phi-weakLone}{\underset{\Leba{\qvd}}{\rightarrow}}\nabt\vw}
    d\vx
    ~dt
    \\\overset{\delta\searrow 0}{\rightarrow}&
    \Big(
    \int_{\R^2} 
    \vv
    \cdot
    \vw
    ~d\vx
    \Big)_{|t=0}
    +
    \int_0^\infty 
    \int_{\R^2} 
    \vv
    \cdot
    \pt\vw
    +
    \vv^2:
    \nabt\vw
    ~d\vx
    ~dt
    \myeqlab{eq:vv-weaksol}
\end{alignat*}

\section{Initial data variety}
\mylabel{section:ini-data}%

\subsection{Locally surjective}

We analyze how $\Vort$ determines $\vrdini$, showing that any $\vrdini\approx\vrdiniz$ can be attained.

\begin{alignat*}{5} \liacs\vrd
    &\topref{eq:vrd}{\eq}
    \liac{\svorf\srad^{\frac1\mv}\sVort} 
    \\&\eq
    \big(\subeq{\liac{\svorf}}{\topref{eq:svorf-lin}{\eq}-\frac1{2\mv}\mdU\lsstf\quad}
    \subeq{\sradz^{\frac1\mv}}{\topref{eq:sradz}{\eq}\mv^{-\frac1{2\mv}}}
    + \subeq{\zeroacc{\svorf}}{\topref{eq:svorfz}{\eq}1}
    \frac1\mv\subeq{\sradz^{\frac1\mv-1} }{\topref{eq:sradz}{\eq}(\mv^{-\frac12})^{\frac1\mv-1}\quad}
    \subeq{\liac{\srad}}{\topref{eq:srad-lin}{\eq}\frac12\mv^{-\frac12}\mdta\lsstf}\big)
    \subeq{\sVortz}{\topref{eq:Vortz}{\eq}\frac{2\mv-1}{\mv}}
    +
    \subeq{\zeroacc{\svorf}}{\topref{eq:svorfz}{\eq}1\quad}
    \subeq{\sradz^{\frac1\mv}}{\topref{eq:sradz}{\eq}(\mv^{-\frac12})^{\frac1\mv}}
    \sVortl
    \\&\eq
    \big(-\frac1{2\mv}\mdU\lsstf~\mv^{-\frac1{2\mv}}
    + \frac1{2\mv} \mv^{-\frac1{2\mv}} \mdta\lsstf \big) \frac{2\mv-1}{\mv}
    + \mv^{-\frac1{2\mv}} \sVortl
    \\&\eq
    \mv^{-\frac1{2\mv}} \Big( \frac{2\mv-1}{2\mv^2}  (\mdta-\mdU)\lsstf + \sVortl \Big)
    \\&\eq
    \mv^{-\frac1{2\mv}} 
    \Big(
    \frac{2\mv-1}{2\mv^2} (\mdta-\mdU)
    \frac{\partial\sstf}{\partial\sVort}(\sstfz,\sVortz)
    \sVortl
    + \sVortl
    \Big)
\intertext{(use the implicit function theorem for Banach spaces)}
    &\eq
    \mv^{-\frac1{2\mv}} 
    \Big( \frac{2\mv-1}{2\mv^2}  (\mdta-\mdU) 
    \big[
    -
    \subeq{\frac{\partial\Fmap}{\partial\sstf}(\sstfz,\sVortz)^{-1}}{\topref{eq:linearization}{\eq}\Lop^{-1}}
    \subeq{\frac{\partial\Fmap}{\partial\sVort}(\sstfz,\sVortz)}{\topref{eq:linearization}{\eq}\moacs2\moacs\mv^2} 
    \big]
    + \moacs1 \Big) \sVortl
    \\&\eq
    \mv^{-\frac1{2\mv}} 
    \Big( (2\mv-1)  (\mdU-\mdta) 
    \Lop^{-1}
    + \moacs1 \Big) \sVortl
    \myeqlab{eq:vort-rad-lin}
\end{alignat*}
Fourier transform of the $T_n$ part: $\sVortl$ is constant in $\ta$, so the transform contains a delta function $\delta(\etav)$:
\begin{alignat*}{5} 
    (\liacs\vrdini)^\wedge(\nv)
    &\topref{eq:vrd-lim-infinity}{=}
    (\liacs\vrd(\ta=0,\uan=\cdot))^\wedge(\nv)
    \\&=
    \int_\R 
    \liacs\vrd^\wedge(\etav,\nv) 
    d\etav
    \\&\toprefc{eq:vort-rad-lin}{eq:mdU-fou}{eq:mdta-fou-a}{\eq}
    \mv^{-\frac1{2\mv}} 
    \int_\R 
    \big(
    \subeq{
        (2\mv-1)  (\Ao+\Bo+\moacs4\moacs\mv) 
        \Lopn^{-1}
    }{\normM{\cdot}\toprefc{eq:I}{eq:B}{eq:A-AB-BA}{\lesssim}\asy\nv^{-1}}
    +\moacs1
    \big)
    \delta(\etav)
    d\etav
    ~(\sVortl)^\wedge(\nv)
    \\&=
    \big(\mv^{-\frac1{2\mv}}+o_\nv(1)\big)
    (\sVortl)^\wedge(\nv)
    \myeqlab{eq:vrdini-tf}
\end{alignat*} 
For $\Nv$ sufficiently large, $(\sVortl)^\wedge\mapsto(\liacs\vrdini)^\wedge$ is close to multiplication by $\mv^{-\frac1{2\mv}}$. 
Now consider $\nv=0$:
\begin{alignat*}{5} \Bo\delta &\topref{eq:mdta-fou-a}{=} \etao\peta\delta = \peta\subeq{\etao\delta}{=0}-\delta \\
    \Ao\delta &\topref{eq:mdU-fou}{=} (\etao-\nvo)\peta \overset{\nv=0}{=} \etao\peta\delta = -\delta \\
    (\Aopn)\delta 
    &\toprefb{eq:mdU-fou}{eq:lstff-new}{=} 
    [(\etao-\nvo)\peta+(\moacs2\pm\moacs\nv)\moacs\mv]\delta 
    \topref{eq:lstff-new}{\underset{\nv=0}{=}} 
    (\etao\peta+\moacs2\moacs\mv)\delta 
    = 
    (2\mv-1)\delta \\
    \Eop\delta &\topref{eq:lstff-new}{=} (2\mv-1)(\Ao-\Bo)\delta = 0 \\
    \Rop\delta &\topref{eq:lstff-new}{=} (\Aop)(\Aon)\Bo\delta = -(2\mv-1)^2\delta \\
    \Lopn\delta &\topref{eq:lstff-new}{=} (\Rop-\Eop)\delta = -(2\mv-1)^2\delta \\ 
    \Lopn^{-1}\delta &= 
    -(2\mv-1)^{-2}\delta \\
    \big((2\mv-1)(\Ao+\Bo+\moacs4\moacs\mv)\Lopn^{-1}+\moacs1\big)\delta 
    &=
    \big((2\mv-1)[(-1)+(-1)+4\mv][-(2\mv-1)^{-2}]+1\big)\delta 
    =
    -\delta
\end{alignat*} 
Hence
\begin{alignat*}{5} 
    \liacs{\vrdini}^\wedge(0) 
    &\topref{eq:vrdini-tf}{=} 
    \mv^{-\frac1{2\mv}} 
    \int_\R [(2\mv-1)(\Ao+\Bo+\moacs4\moacs\mv)\Lopn^{-1}+\moacs1]\delta(\etav) d\etav~\sVortl^\wedge(0)
    = -\mv^{-\frac1{2\mv}} \sVortl^\wedge(0)
    \myeqlab{eq:vrdinizero}
\end{alignat*} 
Again the multiplier $-\mv^{-\frac1{2\mv}}$ is nonzero. Combining the $|\nv|\gg1$ and $\nv=0$ cases we see that $\sVortl\mapsto\liacs\vrdini$ 
is a linear isomorphism on $\Vortdom$. Hence we can take $\frad>0$ sufficiently small to obtain that
$\sVort\mapsto\vrdini$ is a local diffeomorphism on $\Vortnbh$ onto some $\Vortdom$-neighbourhood of $\vrdiniz$.

\subsection{Vortex strength scaling}

\newcommand{\fst}{s}
\newcommand{\fsu}{(\fst)}
\newcommand{\fstf}{\stf^{\fsu}}
\newcommand{\fvv}{\vv^{\fsu}}
\newcommand{\fvort}{\vort^{\fsu}}
\newcommand{\fvrd}{\vrd^{\fsu}}
We have found solutions for initial data $\vrdini$ in some $\Vortdom$-neighbourhood of the constant $\vrdiniz$. 
Neighbourhoods of any other nonzero constant $\vrdinic$ are covered by a simple scaling argument:
if $\vv$ is a weak solution of the incompressible Euler equations, then so is $\fvv(\vx,t):=\fst\vv(\vx,\fst t)$ for any $\fst\in\R$. 
Let $\vv$ be one of the solutions we have already constructed. 
Then 
\begin{alignat*}{5} \fvort(\vx,t):=\nabla_\vx\times\fvv(\vx,t)=\nabla_\vx\times[\fst\vv(\vx,\fst t)]=\fst\vort(\vx,\fst t) \end{alignat*} 
so 
\begin{alignat*}{5} 
    \fvort(\vx,\tv) \rad^{\frac1\mv} 
    \eq 
    \fst \rad^{\frac1\mv} \vort(\vx,\fst\tv) 
    \topref{eq:vrd}{=}
    \fst \vrd(\vx,\fst\tv)
    \overset{t\searrow 0}{\underset{\text{\myeqref{eq:vrd-lim-infinity}}}{\rightarrow}}
    \fst \vrd(\vx,0) 
    =
    \fst \vrdini(\vx) 
\end{alignat*} 
Clearly, given a solution for initial data $\vrdini$, we immediately obtain solutions for 
$\fst\vrdini$ for any $\fst\in\R$.
Hence there is a $\vrad>0$ so that for any constant $\vrdinic$ and initial data $\vrdini$, we have a solution if 
\begin{alignat*}{5} \norma{ \vrdini-\vrdinic }{\Vortdom} < \vrad |\vrdinic| \quad.\end{alignat*}

\section{Pseudo-streamlines algebraic spirals}

\subsection{Shape}

The pseudo-streamlines (one for each fixed $\uan\in\T$) are the graphs of 
\begin{alignat*}{5} \ta\mapsto\vx=\rad
    \begin{bmatrix}
        \cos\pola \\
        \sin\pola
    \end{bmatrix} 
    \toprefb{eq:nonzero-rad}{eq:pola}{=}
    t^\mv\ta^{-\mv}\subeq{\srad}{\topref{eq:nonzero-srad}{\sim}1}
    \begin{bmatrix}
        \cos(\ta+\uan) \\
        \sin(\ta+\uan)
    \end{bmatrix} 
\end{alignat*}
Clearly they remain algebraic spirals 
The spirals do not intersect themselves or each other: we have already shown that our coordinate transforms are bijective,
so $\vx(\ta,\uan)=\vx(\ta_2,\uan)$ for $\ta\neq\ta_2$ (self-intersection) or $\vx(\ta,\uan)=\vx(\ta_2,\uan_2)$ for
$\uan\not\equiv\uan_2$ mod $2\pi$ (two different pseudo-streamlines intersecting) are not possible.

\subsection{Stratification}

\mylabel{section:stratification}

In this section we demonstrate that our solutions are nontrivial in the sense that there is \defm{stratification}
in the inner part of the spiral:
significant variation of vorticity $\vort$ across the near-circular spiral pseudo-streamlines. Stratification could be observed physically, 
so that the algebraic spirals are not a mere mathematical construct. (Vortex sheets, which are not constructed in this paper,
are the extreme case where vorticity not only varies, but vanishes in between curves of concentration.)

Consider one of our solutions for nonzero and non-constant $\Vort$, i.e.\ any solution other than the trivial background solutions.

Consider $(\polb,\Uan)$ coordinates, as in section \myref{section:coordchange}.
Keep $\polb\in\R$ fixed while increasing $\Uan$ from $\Uan_0$ to $\Uan_0+2\pi$. 
This has the effect of moving across pseudo-streamlines radially (see Figure \myref{fig:polb-Uan})
and to return to the same pseudo-streamline ($\uan=\Uan$, the angle of the pseudo-streamline at spatial infinity,
changes by $2\pi$)
but one turn farther from the center ($\ta=\polb-\Uan$, the angle along that streamline, decreases by $2\pi$).
We analyze the variation of $\svort=\svorf\Vort$ using that $\svorf$ varies little, but $\Vort$ varies strongly:
\begin{alignat*}{5} 
    \pU \tao \svorf
    \topref{eq:pU-def}{=}
    (\mdu-\mdta)\tao\svorf
    =
    (\tao(\mdu-\mdta)-\moacs1)\svorf
    =
    (\tao\pU-\moacs1) \svorf
    \topref{eq:mdU}{=}
    (\mdU-\moacs2\moacs\mv) \svorf
    \toprefb{eq:mdU-svorf-bounded}{eq:nonzero-svorf}{\lesssim}
    1
    \myeqlab{eq:pU-vorf}
\end{alignat*} 
For $\Uan_0\leq\Uan_a\leq\Uan_b<\Uan_0+2\pi$,
\begin{alignat*}{5} 
    |(\tao\svorf)(\polb,\Uan_b)-(\tao\svorf)(\polb,\Uan_a)| 
    =
    \Big| \int_{\Uan_a}^{\Uan_b} \pU\tao\svorf~d\Uan \Big|
    \leq
    \int_{\Uan_a}^{\Uan_b} \subeq{|\pU\tao\svorf|}{\topref{eq:pU-vorf}{\lesssim}1} d\Uan
    \lesssim
    |\Uan_b-\Uan_a|
    \leq
    2\pi
    \lesssim 
    1
    \myeqlab{eq:pUvo}
\end{alignat*} 
so
\begin{alignat*}{5} 
    |t\vort(\polb,\Uan_a) - t\vort(\polb,\Uan_b)|
    &\topref{eq:nonzero-vort}{=}
    |(\tao\svorf)(\polb,\Uan_a)\Vort(\Uan_a) - (\tao\svorf)(\polb,\Uan_b)\Vort(\Uan_b)|
    \\&=
    \big|(\tao\svorf)(\polb,\Uan_a)[\Vort(\Uan_a)-\Vort(\Uan_b)] + [(\tao\svorf)(\polb,\Uan_a)-(\tao\svorf)(\polb,\Uan_b)]\Vort(\Uan_b)\big|
    \\&\geq
    \subeq{|(\tao\svorf)(\polb,\Uan_a)|}{\topref{eq:nonzero-svorf}{\sim}\ta}\big|\Vort(\Uan_a)-\Vort(\Uan_b)\big| - \subeq{\big|(\tao\svorf)(\polb,\Uan_a)-(\tao\svorf)(\polb,\Uan_b)\big|}{\topref{eq:pUvo}{\lesssim}1}\subeq{|\Vort(\Uan_b)|}{\topref{eq:nonzero-Vort}{\sim}1}
    \\&\gtrsim
    \ta\big|\Vort(\Uan_a)-\Vort(\Uan_b)\big| - 1 
    \myeqlab{eq:diff-vorf}
\end{alignat*} 
On the other hand 
\begin{alignat*}{5} \tv\vort \topref{eq:nonzero-vort}{\sim} \ta \myeqlab{eq:vortscale} \end{alignat*} 
Taking $\sup$ over $\Uan_a,\Uan_b$ we get
\begin{alignat*}{5} 
    \frac{
        \sup_{\Uan_0\leq\Uan_a\leq\Uan_b<\Uan_0+2\pi}|\vort(\polb,\Uan_a)-\vort(\polb,\Uan_b)|
    }{
        \sup_{\Uan_0\leq\Uan<\Uan_0+2\pi} |\vort(\polb,\Uan)|
    }
    \toprefb{eq:diff-vorf}{eq:vortscale}{\gtrsim}
    \sup_{\Uan_0\leq\Uan_a<\Uan_b<\Uan_0+2\pi}|\Vort(\Uan_a)-\Vort(\Uan_b)| - \ta^{-1}
\end{alignat*} 
which becomes uniformly positive as $\ta\rightarrow\infty$, i.e.\ as we approach the spiral center.

Hence the oscillation of $\vort$ is comparable to or bigger than the size of $\vort$ on each such radial line segment.
This justifies speaking of ``stratified flow'', with vorticity variations that remain observable up to the spiral center.

The proof of Theorem \myref{th:mainresult} and the following remarks is complete.

\nocite{lebeau}
\nocite{caflisch-orellana}
\nocite{lopes-lowengrub-lopes-zheng}
\nocite{duchon-robert}
\nocite{majda-bertozzi-book}
\nocite{delort}
\nocite{yudovich-1995}
\nocite{krasny-icm}
\nocite{sulem-sulem-frisch-bardos}

\bibliographystyle{plain}
\bibliography{../../research/pmeyer/elling}

\end{document}